\documentclass{article}
\usepackage[utf8]{inputenc} % set input encoding (not needed with XeLaTeX)
\usepackage[margin=1.25in]{geometry} % to change the page dimensions
\usepackage{graphicx} % support the \includegraphics command and options
\usepackage{nicefrac}
\usepackage{anyfontsize}

\usepackage{booktabs} % for much better looking tables	
\usepackage{array} % for better arrays (eg matrices) in maths
\usepackage{paralist} % very flexible & customisable lists (eg. enumerate/itemize, etc.)
\usepackage{verbatim} % adds environment for commenting out blocks of text & for better verbatim
\usepackage{mathrsfs}
\usepackage{amsfonts,amsmath,amssymb,amsthm}
\usepackage{esint}
\usepackage{graphics}
\usepackage{enumerate}
\usepackage{mathtools}
\usepackage{xfrac}
\usepackage{bbm}
\usepackage{subcaption}
\usepackage{stmaryrd}

\usepackage[usenames,dvipsnames]{xcolor}
\usepackage[colorlinks=true, pdfstartview=FitV, linkcolor=blue, citecolor=blue, urlcolor=blue]{hyperref}
\usepackage[normalem]{ulem}

\usepackage{tikz}
\usetikzlibrary{calc}
\usepackage{pgf}
\usetikzlibrary{external}
\numberwithin{equation}{section}
\numberwithin{figure}{section}

\newtheorem{theorem}{Theorem}[section]
\newtheorem{lemma}{Lemma}[section]

\theoremstyle{definition}
\newtheorem{definition}[lemma]{Definition}

\newtheorem{remark}[lemma]{Remark}

\newcommand{\ep}{\varepsilon}

\let\tmpsum\sum
\renewcommand{\sum}{\tmpsum\nolimits}
\let\tmplimsup\limsup
\renewcommand{\limsup}{\tmplimsup\nolimits}
\let\tmpsup\sup
\renewcommand{\sup}{\tmpsup\nolimits}

\usepackage{titlesec}

\newcommand{\addperiod}[1]{#1.}
\titleformat{\section}
   {\centering\normalfont\Large}{\thesection.}{0.5em}{}
\titleformat*{\subsection}{\bfseries}
\titleformat{\subsubsection}[runin]
  {\normalfont\bfseries}
  {\thesubsubsection.}
  {0.5em}
  {\addperiod}
\titleformat*{\subsubsection}{\normalfont\itshape}
\titleformat*{\paragraph}{\bfseries}
\titleformat*{\subparagraph}{\large\bfseries}

\title{Stability of a stationary solution to a 1-D model for the MHD}

\author{ % INSERT AUTHOR NAME AFFILIATION AND EMAIL
Yunhao Sun
\thanks{Courant Institute of Mathematical Sciences, New York University.
{\footnotesize \href{mailto: ys5394@nyu.edu}{ys5394@nyu.edu}.}
}
}

\usepackage[nottoc,notlot,notlof]{tocbibind}

\usepackage{enumitem}

\begin{document}

\maketitle

\begin{abstract}
    We investigate the stability of a one-dimensional magnetohydrodynamics model (1-D MHD) with mixed vortex stretching effects, introduced by Dai, Vyas, and Zhang. Using techniques similar to those developed by Lei, Liu, and Ren for the De Gregorio equation, we establish global-in-time well-posedness for initial data near a stationary point. Our result is analogous to the exponential stability of the ground state of the De Gregorio equation.
\end{abstract}

%%%%%%%%%%%%%%%%%%%%%%%%%%%%%%%%%%%%%%%%%%%%%%%%%%%%%%
%                   6. BODY
%%%%%%%%%%%%%%%%%%%%%%%%%%%%%%%%%%%%%%%%%%%%%%%%%%%%%%

% Only the first word and proper nouns of section titles should be capitalized.
% The title of section 1:
\section{Introduction}
%%% INSERT BODY OF PAPER HERE
    In this paper, we will consider the following 1-D MHD model on the torus $\mathbb{T}:=(-\pi, \pi]$
    \begin{align}
        \left\{
        \begin{aligned}
            \omega^+_t + au^-\omega^+_\theta &= p\omega^+H\omega^- + q\omega^-H\omega^+,\\
            \omega^-_t + au^+\omega^-_\theta &= p\omega^-H\omega^+ + q\omega^+H\omega^-,
        \end{aligned}
        \right.\label{1dMHD}
    \end{align}
    with $a, p, q\in \mathbb{R}$ and $H$ being the Hilbert transformation on $\mathbb{T}$ defined by
    \begin{align*}
        Hf(\theta) := \frac{1}{2\pi}\mathrm{p.v.}\int_{-\pi}^\pi\cot\left(\frac{\theta-\vartheta}{2}\right)f(\vartheta)d\vartheta,
    \end{align*}
    while $u^\pm$ are defined to be such that $u^\pm_\theta = H\omega^\pm$ with $u^\pm(0,t)\equiv 0$.
    This model was originally proposed, with scaling $\frac{a}{2} = p = q$, to give a 1-D counterpart to the inviscid 3-D magnetohydrodynamics equations, where $\omega^\pm$ correspond to the Els\"{a}sser variables by mixing vorticity and current \cite{Dai2023}. For a detailed derivation of the 1-D system from the 3-D system, we refer the readers to \cite{Dai2023}. The 1-D MHD is closely related to a number of 1-D models for the incompressible Euler equations, which we will review in the following section.
    \subsection{1-D models for the Euler equations}
    To gain insights to the mechanisms of singularity formation of the 3-D Euler equations, various analogous 1-D models have been proposed and investigated, the first notable one among which being the Constantin-Lax-Majda (CLM) model \cite{Constantin1985}, defined on either $\mathbb{T}$ or $\mathbb{R}$ by
    \begin{equation}
        \omega_t = \omega H\omega.\label{CLM}
    \end{equation}
    In the case of the equation being defined on $\mathbb{R}$, the Hilbert transform is given by
    \begin{align*}
        Hf(x) = \frac{1}{\pi}\mathrm{p.v.}\int_{-\infty}^{+\infty} \frac{f(y)}{x-y} dy.
    \end{align*}
    In the CLM model, $\omega$ mimics the 3-D vorticity while the Hilbert transform of $\omega$ is chosen as a 1-D analog to the 3-D deformation tensor, which gives the effect of vortex stretching. As displayed in \cite{Constantin1985}, equation \eqref{CLM} is equivalent to $\tilde\omega_t = -(i/2) {\tilde\omega}^2$ where $\tilde\omega$ is the holomorphic extension of $\omega$ defined by $\tilde{\omega} := \omega + iH\omega$ onto either the unit disk or the half plane. Hence, local well-posedness is guaranteed for initial data satisfying $\lVert \omega_0 \rVert_{L^\infty} + \lVert H\omega_0 \rVert_{L^\infty} < \infty$ and finite time blowup occurs if and only if the initial datum $\omega_0 + iH\omega_0$ samples purely imaginary values.
    \par
    Notice that in the CLM model, only the vortex stretching effect is present, this was pointed out in \cite{DeGregorio1990}. To compare the effects of advection and vortex stretching, the generalized Constantin-Lax-Majda model (gCLM) was proposed \cite{Wunsch2008}, which is defined (on either $\mathbb{T}$ or $\mathbb{R}$) by
    \begin{equation}
        \omega_t + au \omega_x = \omega H\omega,\quad\text{where}\quad u_x = H\omega,
    \end{equation}
    where $a\in \mathbb{R}$ is some fixed constant. Singularity formation of gCML has been extensively researched in literature. In the case of $a<0$, studies in \cite{Castro2010} showed that advection and vortex stretching work together in producing finite in time blowup of classical solutions on $\mathbb{R}$ given certain $C_c^\infty$ data. The case of $a=0$ reduces to CLM. While for small $a>0$, results in \cite{Elgindi2020} showed that there exist $H^3$ initial data on $\mathbb{R}$ that produce singularity in finite time, and moreover, H\"{o}lder initial data on $\mathbb{R}$ which produce singularity were found for any given $a\in\mathbb{R}$. The results for $a>0$ were further explored in \cite{Chen2021}, where blowup was proven for certain $C^\infty_c$ initial data on $\mathbb{R}$ given $a=1$ as well as $C^\infty_c$ data on $\mathbb{T}$ given small $a>0$.
    \par
    For the $a=1$ case, the gCLM model is known as the De Gregorio equation \cite{DeGregorio1990}. In comparison to singularity formation, global well-posedness of the De Gregorio equation also attracted vast attention. In \cite{Lei2020}, the authors established global in time well-posedness on both $\mathbb{T}$ or $\mathbb{R}$ with positive compact supported initial data in $H^k$ for $k\geq1$. In \cite{Chen2023}, the author studied the De Gregorio equation on $\mathbb{T}$ with odd $H^1$ initial data nonnegative on $[0,\pi]$ and discovered a one-point continuation criterion as well as a subset of such initial data that admits global existence. Furthermore, global in time solutions on $\mathbb{T}$ that converge to stationary points are also shown to exist. We will review these results in the next section.
    \subsection{Convergence of the De Gregorio equation to the equilibrium \texorpdfstring{$-\sin(\theta)$}{-sin(theta)}}
    Recall the De Gregorio model on the torus defined by
    \begin{equation}
        \omega_t + u \omega_\theta = \omega H\omega,\quad\text{where}\quad u_\theta = H\omega\quad\text{and}\quad u(0, t) \equiv 0,\label{DeGregorio}
    \end{equation}
    where we have chosen a gauge of $u$ by fixing $u(0,t)\equiv 0$. It is evident that $\omega = \Omega := -\sin(\theta)$ with $u = U := \sin(\theta)$ gives a time-stationary solution to \eqref{DeGregorio}. In \cite{Jia2019}, the authors considered perturbing the stationary solution by letting $\omega = \Omega + \ep\eta$ and showed that given mean zero initial perturbation $\eta_0$ in $H^\frac{3}{2}$ and the weighted $L^2$ space $\tilde Y$ defined by
    \begin{align*}
        \tilde Y := \left\{f\in L^2(\mathbb{T}) : \int_\mathbb{T} \frac{|f|^2}{|{\sin(\theta/2)}|^{2\gamma}}d\theta < \infty\right\},
    \end{align*}
    where $\gamma \in (3/2, 2)$, for $\ep$ small enough, the $H^\frac{3}{2}$ norm of the perturbed term $\eta$ remains bounded in time, and $\eta$ converges to zero exponentially in $\tilde Y$ and $H^s$ for any $s<3/2$. Furthermore, in \cite{Lei2020}, the authors proposed an alternative weighted Sobolev space $\mathcal{H}$ defined by
    \begin{align}
        \mathcal{H} := \left\{f\in H^1(\mathbb{T}) : f(0) = 0\quad\text{and}\quad\int_\mathbb{T} \frac{|f_\theta|^2} {|
        {\sin(\theta/2)}|^2}d\theta < \infty\right\}
    \end{align}
    and showed that for initial perturbation $\eta_0 \in \mathcal{H}$ with small $\ep$, the solution $\omega = \Omega + \ep\eta$ converges exponentially in $\mathcal{H}$ to a point that depends on the initial average of $\eta_0$.
    \par
    Observe that the 1-D MHD \eqref{1dMHD} has the property of when taken initial data $\omega^+_0=\omega^-_0$, the model reduces to the gCLM model with $\omega = \omega^+ = \omega^-$, or the De Gregorio equation when further taken $a = p+q$. This suggests convergence results similar to the above might also hold true for the 1-D MHD. In the next section, we will recall some fundamental properties of the 1-D MHD discussed in \cite{Dai2023}.
    \subsection{Local well-posedness and BKM type criterion of the 1-D MHD}
    In \cite{Dai2023}, the authors have established $H^1$ local in time well-posedness in the case of $p = 0, q = 1$ and arbitrary $a \in \mathbb{R}$ by appealing to the Kato-Lai existence theorem \cite{Kato1984}. They also found a BKM type of continuation criterion by controlling $\lVert H\omega^+ \rVert_{L^\infty} + \lVert H\omega^- \rVert_{L^\infty}$. Finally, they have also demonstrated global well-posedness when there are no vortex stretching terms, i.e., $p = q = 0$. We will illustrate an alternative proof for local existence with arbitrary $a$, $p$, and $q$, as well as a continuation criterion, outlined in \ref{existence}.
    \subsection{Reformulation of the 1-D MHD perturbed at the equilibrium \texorpdfstring{$-\sin(\theta)$}{-sin(theta)}}\label{reformulation}
    In this paper, we will consider \eqref{1dMHD} with $a = 1$ and $p + q = 1$ with $p$, $q$ being positive, the scaling of these parameters gives positive vortex stretching terms and preserves the Lie bracket structure of the MHD, which can be observed in \eqref{Leta+} and \eqref{Leta-}. In considering the behavior of $\omega^\pm$ near the stationary point $\Omega=-\sin(\theta)$, we will rewrite \eqref{1dMHD} in terms the variables $\eta^\pm$ defined on $\mathbb{T}$ by
    \begin{align*}
        \omega^\pm(\theta,t) = -\sin(\theta)+\ep (\eta^+(\theta, t)\pm\eta^-(\theta, t)),
    \end{align*}
    we also define accordingly $v^\pm$ as such that $\partial_\theta v(\theta, t) = H\eta^\pm(\theta, t)$ with $v^\pm(0, t)\equiv 0$. Hence, equation \eqref{1dMHD} with $\omega^\pm$ defined in terms of $\eta^\pm$ gives the following PDE of $\eta^\pm$
    \begin{align}
        \eta^+_t &= \{\eta^++v^+, \sin(\theta)\} + \ep N_1,\label{Leta+}\\
        \eta^-_t &= \{\eta^--v^-, \sin(\theta)\} - 2q(\cos(\theta)\eta^- + \sin(\theta)H\eta^-) + \ep N_2,\label{Leta-}\\
        N_1 &= \{\eta^+, v^+\} - \{\eta^-, v^-\},\notag\\
        N_2 &= \{\eta^-, v^+\} - \{\eta^+, v^-\} - 2q(\eta^-H\eta^+ - \eta^+H\eta^-)\notag,
    \end{align}
    where $\{\cdot, \cdot\}$ denotes the Lie bracket, i.e. $\{f, g\} = fg_\theta - f_\theta g$, and $N_1, N_2$ denote the nonlinear terms. It is worth noting that the linear parts of the coupled PDE are decoupled, and in comparison to the perturbed De Gregorio model around $\Omega = -\sin(\theta)$, which is given in \cite{Jia2019, Lei2020} by
    \begin{align*}
        \eta_t = \{\eta + v, \sin(\theta)\} + \ep\{\eta, v\},
    \end{align*}
    the linear operator on $\eta^+$ is identical to the above, while for $q=0$ the linear operator on $\eta^-$ differs from the above by a multiple of the term $\{v^-,\sin(\theta)\}$. It is not hard to see that this difference is a bounded operator on $\eta^-$ with respect to $L^2(\mathbb{T})$ and any Sobolev norm $H^s(\mathbb{T})$ for $s>1/2$. These observations motivate theorem \ref{thm1}.
    \par
    In the remaining passage, we will use the following definitions on the linear operators.
    \begin{definition}
        Define $L$, $B$, and $Q$ to be the operators
        \begin{equation}
            Lf = \{f,\sin(\theta)\},\quad Bf = \{v(f), \sin(\theta)\},\quad\text{and}\quad Qf = \cos(\theta)f+\sin(\theta)Hf, \label{definition}
        \end{equation}
        where $v(f)$ denote the function $v$ on $\mathbb{T}$ such that $v_\theta = Hf$ with $v(0) = 0$. Also, let $L^+ := L+B$ and $L^-:=L-B-2qQ$, then equation \eqref{Leta+} and \eqref{Leta-} can be written as
        \begin{align*}
            \eta^+_t &= L^+\eta^+ + \ep N_1 = (L+B)\eta^+ + \ep N_1,\\
            \eta^-_t &= L^-\eta^- + \ep N_2 = (L-B-2qQ)\eta^- + \ep N_2.
        \end{align*}
    \end{definition}
    \subsection{Main results}
    \begin{theorem}\label{thm1}
        Consider the 1-D MHD in the case $a=1$, $p = 1$, $q=0$, the stable in time solution $\Omega = -\sin(\theta)$ and $U = \sin(\theta)$, and the weighted Sobolev space $\mathcal{H}$ defined above. Given any mean zero $H^2$ initial data $\omega^\pm_0$ close enough to $\Omega$ in $\mathcal{H}$, i.e., with $\omega^\pm_0 - \Omega$ small in the $\mathcal{H}$ norm, we have that $\omega^\pm(\cdot, t)$ converge in time exponentially to $\Omega$ in the $\mathcal{H}$ norm.
    \end{theorem}
    \begin{theorem}\label{thm2}
        In the case of $a=1$, $p+q=1$, $0<q<\frac{1}{4}$. For any mean zero $H^2$ initial data $\omega^\pm_0$ close enough to $\Omega$ in $\mathcal{H}$, there exists a bounded continuous function $h:\mathbb{R}_+\rightarrow \mathbb{R}$ such that 
        the differences between $\omega^\pm(\cdot, t)$ and $(1 \pm h)\Omega$ decay exponentially in $\mathcal{H}$ in time. Furthermore, $h(t)$ converges as $t\rightarrow\infty$ and $\sup_{t \geq 0}|h(t)|\in\mathcal{O}(\ep q)$ for $\ep, q\rightarrow 0$.
    \end{theorem}
    \begin{remark}
        Our proof idea is straightforward, for the linear operators, in section \ref{linearq=0}, similar to the approach of \cite{Jia2019} and \cite{Lei2020}, we study the anti-Hermitian property of the operators $L^\pm$ with respect to a basis of $\mathcal{H}$ and establish decay estimates. For the nonlinear components, in section \ref{nonlinearq=0}, we use joint analysis techniques on $N_1$ and $N_2$ analogous to those for the 3-D magnetohydrodynamics equations. For the case of $0<q<\frac{1}{4}$, we show that the decay effect of $L^\pm$ dominates the effect of the operator $Q$, with details in section \ref{whyq<1/4}.
    \end{remark}
    \section{An anisotropic Sobolev space}
    We will analyze the perturbed 1-D MHD in the following function space inspired by \cite{Jia2019, Lei2020}
    \begin{align*}
        Y:=\operatorname{span}_{\mathbb{R}}\{\sin(\theta)\} \oplus \left\{f\in H^1(\mathbb{T}) : f(0) = 0, \int_{\mathbb{T}} \frac{|f_\theta|^2}{\left|{\sin(\theta/2)}\right|^2}d\theta < \infty\right\}.
    \end{align*}
    Denote the spaces above respectively by $Z_1$ and $\mathcal{H}$ so $Y = Z_1\oplus \mathcal{H}$. It is shown in \cite{Lei2020} that the following functions $e_{c,0}$, $e_{s,k}$'s, and $e_{c,k}$'s for $k\in\mathbb{Z}_+$ on $\mathbb{T}$ defined by $e_{c, 0} := \cos(\theta) - 1$ and
    \begin{align*}
        \begin{aligned}
            e_{s, k} &:= \frac{\sin\left((k+1)\theta\right)}{k+1} - \frac{\sin (k\theta)}{k} &&\quad\text{for}\quad k\geq 1,\\
            e_{c, k} &:= \frac{\cos\left((k+1)\theta\right) - 1}{k+1} - \frac{\cos (k\theta) - 1}{k} &&\quad\text{for}\quad k\geq 1,
        \end{aligned}
    \end{align*}
    gives an orthonormal basis of $\mathcal{H}$. We will decompose $\mathcal{H}$ into two subspaces,
    \begin{align*}
        Z_2:=\operatorname{span}_\mathbb{R}\{\cos(\theta)-1\} \quad\text{and}\quad \mathcal{H}_0:= \overline{\operatorname{span}}_\mathbb{R}\{e_{s,k} : k\geq 1\} \oplus \overline{\operatorname{span}}_\mathbb{R}\{e_{c,k} : k\geq 1\}.
    \end{align*}
    Then given any function $f$ in $Y$, we can write $f$ as a sum with coefficients $f_{s,k}$'s and $f_{c,k}$'s
    \begin{align*}
        f(\theta) = f_{s,0}\sin(\theta) + f_{c,0}\left(\cos(\theta)-1\right) + \sum\nolimits_{k\geq 1}\left(f_{s,k}e_{s,k}(\theta) + f_{c,k}e_{c,k}(\theta)\right),
    \end{align*}
    and for $f, g\in Y$, we define the inner product on $Y$ by
    \begin{align*}
        \langle f, g\rangle_Y := f_{s,0}g_{s,0} + \langle \mathbb{P}_\mathcal{H} f, \mathbb{P}_\mathcal{H} g\rangle_\mathcal{H} := f_{s,0}g_{s,0} + \frac{1}{4\pi}\int_{\mathbb{T}} \frac{\partial_\theta(\mathbb{P}_{\mathcal{H}}f) \partial_\theta(\mathbb{P}_{\mathcal{H}}g)}{\left|{\sin(\theta/2)}\right|^2} d\theta,
    \end{align*}
    where $\mathbb{P}_\mathcal{H}$ is the projection operator onto $\mathcal{H}$, and for $f\in Y$, the norm is defined accordingly by
    \begin{align*}
        \lVert f \rVert_Y^2 := \llbracket f \rrbracket_{Z_1}^2 + \llbracket f \rrbracket_{\mathcal{H}}^2 = |f_{s,0}|^2 + \frac{1}{4\pi}\int_{\mathbb{T}} \left|\frac{\partial_\theta(\mathbb{P}_{\mathcal{H}}f)}{\sin\left(\theta/2\right)}\right|^2 d\theta.
    \end{align*}
    We will also decompose the $\mathcal{H}$ seminorm into $\llbracket \cdot \rrbracket_{Z_2}$ and $\llbracket \cdot \rrbracket_{\mathcal{H}_0}$ given by
    \begin{align*}
        \llbracket f \rrbracket_{Z_2}^2 := |f_{c,0}|^2 = \langle f, \mathbb{P}_{Z_2} f\rangle_{\mathcal{H}}\quad\text{and}\quad\llbracket f \rrbracket_{\mathcal{H}_0}^2 = \langle f, f\rangle_{\mathcal{H}_0} := \langle f, \mathbb{P}_{\mathcal{H}_0} f\rangle_{\mathcal{H}},
    \end{align*}
    where we have denoted $\langle\cdot, \cdot\rangle_{\mathcal{H}_0}$ to be the inner product on $\mathcal{H}_0$ induced from $\mathcal{H}$. With a slight abuse of notation, we will also often consider $\langle\cdot, \cdot\rangle_{\mathcal{H}}$ and $\langle\cdot, \cdot\rangle_{\mathcal{H}_0}$ as pseudo inner products on $Y$.
    \subsection{Properties of the PDE operators with respect to the subspaces of \texorpdfstring{$Y$}{Y}}
    Here, we will point out a few observations of the operators $L$ and $B$ with respect to $Y$. It is clear that $B$ is a bounded operator on $Y$.
    $L\pm B$ are closed and densely defined on $Y$. $L$ and $B$ map $Z_1$ to zero, and $Z_2$ is invariant under $L$ and $B$, these claims can be observed from the section on linear analysis that follows.
    \subsection{Some conserved quantities}
    \subsubsection{Conservation of the average of \texorpdfstring{$\eta^+$}{eta+}} \label{conserve}
    By a direct integration of $\eta^+$ on $\mathbb{T}$, we see that
    \begin{align*}
        \int_{\mathbb{T}}\eta^+(\theta, t)d\theta \equiv \int_{\mathbb{T}}\eta^+(\theta, 0)d\theta\quad\text{for any}\quad t\geq0.
    \end{align*}
    Observe that the averages of basis functions $e_{s,k}$'s are zero but the averages of $e_{c,k}$'s are nonzero, hence we can get the following relation between $\eta^+_{c,k}$'s
    \begin{align*}
        2\pi\sum\nolimits_{k\geq1}\tfrac{1}{k(k+1)}\eta^+_{c,k} - 2\pi\eta^+_{c,0} = \int_{\mathbb{T}}\eta^+(\theta, 0)d\theta.
    \end{align*}
    Hence if we start with initial data $\eta^+_0$ with zero average, then we have the equivalence
    \begin{align*}
        \llbracket \eta^+ \rrbracket_{\mathcal{H}_0}\leq \llbracket \eta^+ \rrbracket_\mathcal{H} \lesssim \llbracket \eta^+ \rrbracket_{\mathcal{H}_0}.
    \end{align*}
    \subsubsection{Conservation of \texorpdfstring{$\eta^\pm(0, t)$}{eta(0, t)}} Let $Z_0 := \operatorname{span}_\mathbb{R}\{1\}$ be the space of constant, which is the orthogonal complement of $Y$ in $L^2(\mathbb{T})$. Given \ref{embedding}, we see that $Z_0\oplus Y\hookrightarrow H^1(\mathbb{T})\hookrightarrow C(\mathbb{T})$, hence the pointwise values of $\eta^\pm(0,t)$ and $H\eta^\pm (0,t)$ are well-defined for all $\eta^\pm\in Z_0\oplus Y$ and we have
    \begin{align*}
        \frac{d}{dt}\eta^+(0,t) =& \eta^+(0,t) + \ep H\eta^+(0,t)\eta^+(0,t) - \ep H\eta^-(0,t)\eta^-(0,t),\\
        \frac{d}{dt}\eta^-(0,t) =& \eta^-(0,t) - 2q\eta^-(0,t) + \ep H\eta^+(0,t)\eta^-(0,t) - \ep H\eta^-(0,t)\eta^+(0,t)\\
        &- 2q\ep(H\eta^+(0,t)\eta^-(0,t) - H\eta^-(0,t)\eta^+(0,t)).
    \end{align*}
    It is apparent that $\eta^\pm(0,t)$ is equivalent to the $Z_0$ component of $\eta^\pm \in Z_0\oplus Y$, and if we start with $\eta^\pm_0\in Y$, i.e., $\eta^\pm_0(0,0) = 0$, this remains true for all $t\geq 0$. This justifies our choice of $Y$.
    \subsubsection{Conservation of \texorpdfstring{$\eta^\pm_\theta(0,t)$}{eta_theta(0,t)}}
    By the Sobolev embedding $H^2(\mathbb{T})\hookrightarrow C^1(\mathbb{T})$, we see that for all $\eta^\pm\in H^2(\mathbb{T})$, the pointwise values of $\eta^\pm_\theta(0, t)$ and $H\eta^\pm_\theta(0, t)$ are well-defined, and we have
    \begin{align*}
        \frac{d}{dt}\eta^+_\theta(0, t) &= \ep H\eta^+_\theta(0, t)\eta^+(0, t) - \ep H\eta^-_\theta(0,t)\eta^-(0, t)\\
        \frac{d}{dt}\eta^-_\theta(0, t) &= -2q(\eta^-_\theta(0, t) + H\eta^-(0, t)) + \ep H\eta^+_\theta(0, t)\eta^-(0, t) - \ep H\eta^-_\theta(0,t)\eta^+(0, t)\\
        &\quad - 2q\ep (H\eta^+(0, t)\eta^-_\theta(0, t)+H\eta^+_\theta(0, t)\eta^-(0, t)-H\eta^-(0,t)\eta^+_\theta(0,t)\\
        &\quad -H\eta^-_\theta(0,t)\eta^+(0,t)).
    \end{align*}
    For $\eta^\pm\in H^2(\mathbb{T})$, the $Z_0$ component of $\eta^\pm$ corresponds to $\eta^\pm(0, t)$ and the $Z_1$ component corresponds to $\eta^\pm_\theta(0, t)$. Hence for $\eta^\pm_0\in H^2(\mathbb{T})\cap \mathcal{H}$, in the case of $q=0$, the $Z_0$ and $Z_1$ component of $\eta^\pm$ will remain zero. This justifies our use of $\mathcal{H}$ in the case of $q=0$. For $q>0$ however, we have
    \begin{equation}
        \frac{d}{dt}\eta^-_\theta(0, t) = -2q(\eta^-_\theta(0, t) + H\eta^-(0, t)) - 2q\ep H\eta^+(0, t)\eta^-_\theta(0, t).\label{etathetapointwise}
    \end{equation}
    \section{Linear stability in the case of \texorpdfstring{$q=0$}{q=0}}\label{linearq=0}
    In this section, we will primarily focus on the linear estimates of $\eta^\pm$ in the $\mathcal{H}_0$ seminorm in the $q=0$ case with $\eta^\pm_0 \in \mathcal{H}$. To recall the linear part of the perturbed 1-D MHD
    \begin{align}
        \eta^+_t &= L^+\eta^+  = (L+B) \eta^+ = \{\eta^++v^+, \sin(\theta)\},\\
        \eta^-_t &= L^-\eta^- = (L-B) \eta^- = \{\eta^--v^-, \sin(\theta)\}.
    \end{align}
    \subsection{Linear estimate on \texorpdfstring{$\eta^+$}{eta+}}
    As mentioned in section \ref{reformulation}, the linear operator $L^+$ on $\eta^+$ is identical to the linear operator of the perturbed De Gregorio equation. Here we will first summarize the results shown in \cite{Jia2019, Lei2020}. The operator $L^+$ on the Fourier basis defined by trigonometric functions $\{\sin(k\theta) : k\in \mathbb{Z}_+\} \cup \{\cos(k\theta) : k\in \mathbb{Z}_+\} \cup \{1\}$ gives $L^+ 1 = \cos(\theta)$, and for $k\geq 1$
    \begin{align*}
        L^+ \sin (k\theta) &= a^+_k \sin \left((k-1)\theta\right) - b^+_k \sin \left((k+1)\theta\right),\\
        L^+ \cos (k\theta) &= a^+_k\cos \left((k-1)\theta\right) - b^+_k\cos \left((k+1)\theta\right) + \tfrac{1}{k}\cos (\theta),
    \end{align*}
    where $a^+_k$ and $b^+_k$ are coefficients defined by
    \begin{align*}
        a^+_k:= \frac{1}{2}(k+1)\left(1-\tfrac{1}{k}\right),\qquad
        b^+_k := \frac{1}{2}(k-1)\left(1-\tfrac{1}{k}\right).
    \end{align*}
    A direct calculation shows that, on the basis $\{e_{s,k} : k\in\mathbb{N}\}\cup\{ e_{c,k} : k\in \mathbb{N}\}$, $L^+e_{s,0} = L^+ e_{c,0} = 0$, and
    \begin{align*}
        L^+e_{s,k} &= -d^+_{k+1}e_{s,k+1} - (d^+_{k+1}-d^+_k)e_{s,k} + d^+_ke_{s,k-1},\\
        L^+e_{c,k} &= -d^+_{k+1}e_{c,k+1} - (d^+_{k+1}-d^+_k)e_{c,k} + d^+_ke_{c,k-1}+\tfrac{k^2-k-1}{k^2(k+1)^2} e_{c,0},
    \end{align*}
    where $d^+_k=\frac{(k-1)^2(k+1)}{2k^2}$ for $k\geq 1$ and $d^+_{k+1}-d^+_k\geq\frac{3}{8}$ for all $k\geq 1$. Hence for $\eta^+\in Y$, we have
    \begin{align*}
        (L^+\eta^+)_{s, k} &= d^+_{k+1}\eta^+_{s,k+1} - (d^+_{k+1} - d^+_k)\eta^+_{s,k} - d^+_k\eta^+_{s,k-1} \quad\text{for}\quad k\geq 1,\\
        (L^+\eta^+)_{c, k} &= d^+_{k+1}\eta^+_{c,k+1} - (d^+_{k+1} - d^+_k)\eta^+_{c,k} - d^+_k\eta^+_{c,k-1} \quad\text{for}\quad k\geq 1,\\
        (L^+\eta^+)_{c,0} &= \sum\nolimits_{k\geq 1}\tfrac{k^2-k-1}{k^2(k+1)^2} \eta^+_{c,k}.
    \end{align*}
    \subsection{Linear estimate on \texorpdfstring{$\eta^-$}{eta-}}
    For the second equation, for $q=0$, $L^-$ on the Fourier basis gives $L^-1=\cos(\theta)$, and for $k\geq1$
    \begin{align*}
        L^- \sin (k\theta) &= a^-_k \sin \left((k-1)\theta\right) - b^-_k \sin \left((k+1)\theta\right),\\
        L^- \cos (k\theta) &= a^-_k\cos \left((k-1)\theta\right) - b^-_k\cos \left((k+1)\theta\right) - \tfrac{1}{k}\cos (\theta),
    \end{align*}
    where $a^-_k$ and $b^-_k$ are given by    \begin{align*}
        a^-_k:=\frac{1}{2}(k+1)\left(1+\tfrac{1}{k}\right),\qquad
        b^-_k:=\frac{1}{2}(k-1)\left(1+\tfrac{1}{k}\right),
    \end{align*}
    and on the basis $\{e_{s,k} : k\in\mathbb{N}\}\cup\{ e_{c,k} : k\in \mathbb{N}\}$ of $Y$, $L^-e_{s,0} = 0$, $L^- e_{c,0} = -2e_{c,0}$,
    \begin{align*}
        L^-e_{s,k} &= -d^-_{k+1}e_{s,k+1} - (d^-_{k+1}-d^-_k)e_{s,k} + d^-_ke_{s,k-1},\\
        L^-e_{c,k} &= -d^-_{k+1}e_{c,k+1} - (d^-_{k+1}-d^-_k)e_{c,k} + d^-_ke_{c,k-1} +\tfrac{k^2+3k+1}{k^2(k+1)^2}e_{c, 0},
    \end{align*}
    where $d_k =\frac{(k-1)(k+1)^2}{2k^2}$ for $k\geq 1$, and note that the terms $d^-_{k+1} - d^-_k$ satisfy
    \begin{align*}
        d^-_{k+1} - d^-_k = \frac{1}{2} + \tfrac{k^2+3k+1}{2k^2(k+1)^2} \geq \frac{1}{2}\quad\text{for all}\quad k\geq1.
    \end{align*}
    So, for $\eta^-\in Y$, we have that
    \begin{align*}
        (L^-\eta^-)_{s, k} &= d^-_{k+1}\eta^-_{s,k+1} - (d^-_{k+1} - d^-_k)\eta^-_{s,k} - d^-_k\eta^-_{s,k-1} \quad\text{for}\quad k\geq 1,\\
        (L^-\eta^-)_{c, k} &= d^-_{k+1}\eta^-_{c,k+1} - (d^-_{k+1} - d^-_k)\eta^-_{c,k} - d^-_k\eta^-_{c,k-1} \quad\text{for}\quad k\geq 1,\\
        (L^-\eta^-)_{c,0} &= -2\eta^-_{c,0} + \sum\nolimits_{k\geq 1}\tfrac{k^2+3k+1}{k^2(k+1)^2} \eta^-_{c,k}.
    \end{align*}
    Hence for the linear estimates, we have the following conclusions.
    \begin{lemma}
        For $\eta^\pm\in Y$, we have the following linear estimate when $q=0$,
        \begin{align*}
            \langle\eta^+, L^+\eta^+\rangle_{\mathcal{H}_0} \leq -\frac{3}{8}\llbracket \eta^+ \rrbracket_{\mathcal{H}_0}^2 ,\quad \langle\eta^-, L^-\eta^-\rangle_{\mathcal{H}_0} \leq -\frac{1}{2}\llbracket \eta^- \rrbracket_{\mathcal{H}_0}^2
        \end{align*}
        and for the space $Z_2$, the linear evolution gives
        \begin{align*}
            (L^+\eta^+)_{c,0} = \sum\nolimits_{k\geq 1}\tfrac{k^2-k-1}{k^2(k+1)^2} \eta^+_{c,k},\quad (L^-\eta^-)_{c,0} = -2\eta^-_{c,0} + \sum\nolimits_{k\geq 1}\tfrac{k^2+3k+1}{k^2(k+1)^2} \eta^-_{c,k}.
        \end{align*}
        \begin{proof}
            The time derivative of $\eta^\pm_{c,0}$ is already derived above. To show the first claim,
            \begin{align*}
                \langle \eta^+, L^+ \eta^+\rangle_{\mathcal{H}_0} =& \sum\nolimits_{k\geq 1}\eta^+_{s,k}\left(d^+_{k+1}\eta^+_{s,k+1} - (d^+_{k+1} - d^+_k)\eta^+_{s,k} - d^+_k\eta^+_{s,k-1}\right)\\
                &+ \sum\nolimits_{k\geq 1}\eta^+_{c,k}\left(d^+_{k+1}\eta^+_{c,k+1} - (d^+_{k+1} - d^+_k)\eta^+_{c,k} - d^+_k\eta^+_{c,k-1}\right)\\
                =& \sum\nolimits_{k\geq 1}(d^+_{k+1} - d^+_k)\left((\eta^+_{s,k})^2 + (\eta^+_{c,k})^2\right)\\
                \leq& -\frac{3}{8}\llbracket \eta^+ \rrbracket_{\mathcal{H}_0}^2.
            \end{align*}
            Note that the sum above converges absolutely on $D(L)\cap Y$, then we apply the continuity of the quadratic form $\eta^+\mapsto\langle\eta^+, L^+\eta^+\rangle_{\mathcal{H}_0}$. The estimate on $L^-$ holds analogously.
        \end{proof}
    \end{lemma}
    \section{{Nonlinear analysis in the case of \texorpdfstring{$q=0$}{q=0}}}\label{nonlinearq=0}
    In this section, we will consider the nonlinear estimates of the perturbed 1-D MHD in the case of $q=0$ and conclude the proof of \ref{thm1}. First recall the nonlinear perturbed 1-D MHD when $q=0$
    \begin{align}
        \eta^+_t &= L^+\eta^+ + \ep \{\eta^+, v^+\} - \ep\{\eta^-, v^-\},\label{q=0eta+}\\
        \eta^-_t &= L^-\eta^- + \ep \{\eta^-, v^+\} - \ep\{\eta^+, v^-\}.\label{q=0eta-}
    \end{align}
    \subsection{Estimates of \texorpdfstring{$\eta^+$}{eta+} and \texorpdfstring{$\eta^-$}{eta-} in \texorpdfstring{$\llbracket \cdot \rrbracket_{\mathcal{H}_0}$}{H0}}
    To estimate $\eta^+(\cdot, t)$ in $\llbracket \cdot \rrbracket_{\mathcal{H}_0}$ for $\eta^+\in \mathcal{H}$,
    \begin{align*}
        \frac{1}{2}\frac{d}{dt}\llbracket \eta^+ \rrbracket_{\mathcal{H}_0}^2 = \langle \eta^+, L^+\eta^+\rangle_{\mathcal{H}_0} + \ep\langle \{\eta^+, v^+\}, \eta^+\rangle_{\mathcal{H}_0} - \ep\langle\{\eta^-, v^-\}, \eta^+\rangle_{\mathcal{H}_0},
    \end{align*}
    for the nonlinear part, we will write it in terms of three terms
    \begin{align*}
        I^+_1 &:= \frac{1}{4\pi}\int_{\mathbb{T}}\frac{1}{\left|{\sin(\theta/2)}\right|^2}\left(\eta^+v^+_{\theta\theta}\eta^+_\theta - \eta^+_{\theta\theta}v^+\eta^+_\theta - \eta^-v^-_{\theta\theta}\eta^+_\theta \right)d\theta,\\
        I^+_2 &:= \frac{1}{4\pi} \int_{\mathbb{T}}\frac{\eta^+_{c,0}\sin(\theta)}{\left|{\sin(\theta/2)}\right|^2}\left(\eta^+v^+_{\theta\theta} - \eta^+_{\theta\theta}v^+ - \eta^-v^-_{\theta\theta} + \eta^-_{\theta\theta}v^- \right)d\theta,\\
        I^+_3 &:= \frac{1}{4\pi}\int_{\mathbb{T}}\frac{1}{\left|{\sin(\theta/2)}\right|^2} \left(\eta^-_{\theta\theta}v^-\eta^+_\theta \right)d\theta,
    \end{align*}
    It is worth noting that $I_1 + I_3 = \langle N_1, \eta^+\rangle_\mathcal{H}$ and $I_2 = \langle N_1, \mathbb{P}_{Z_2}\eta^+\rangle_\mathcal{H}$, hence $I_1-I_2+I_3 = \langle N_1, \eta^+\rangle_{\mathcal{H}_0}$. To bound $I_1$, note that $\lVert v^\pm_{\theta\theta}\rVert_{L^2} \leq \lVert H\eta^\pm \rVert_{H^1}\leq \llbracket \eta^\pm \rrbracket_{\mathcal{H}}$, then
    \begin{align*}
        |I^+_1| &\lesssim \left\lVert\frac{\eta^+}{{\sin(\theta/2)}}\right\rVert_{L^\infty}\lVert v^+_{\theta\theta }\rVert_{L^2}\llbracket \eta^+ \rrbracket_{\mathcal{H}} + \left\lVert\frac{\eta^-}{{\sin(\theta/2)}}\right\rVert_{L^\infty}\lVert v^-_{\theta\theta }\rVert_{L^2}\llbracket \eta^+ \rrbracket_{\mathcal{H}}\\
        &\quad + \left|\int_{\mathbb{T}}\frac{v^+\partial_\theta(\eta^+_\theta)^2}{\left|{\sin(\theta/2)}\right|^2}d\theta\right|\\
        &\lesssim \llbracket \eta^+ \rrbracket_{\mathcal{H}}^3 + \llbracket \eta^+ \rrbracket_{\mathcal{H}}\llbracket \eta^- \rrbracket_{\mathcal{H}}^2 + \left|\int_{\mathbb{T}}\frac{v^+_\theta\eta^+_\theta\eta^+_\theta}{\left|{\sin(\theta/2)}\right|^2}d\theta\right|
        + \left|\int_{\mathbb{T}} \frac{\cos\left(\theta/2\right)}{{\sin(\theta/2)}^3} \left(v^+\eta^+_\theta\eta^+_\theta\right) d\theta\right|\\
        &\lesssim \llbracket \eta^+ \rrbracket_{\mathcal{H}}^3 + \llbracket \eta^+ \rrbracket_{\mathcal{H}}\llbracket \eta^- \rrbracket_{\mathcal{H}}^2 + \lVert H\eta^+ \rVert_{L^\infty}\llbracket \eta^+ \rrbracket_{\mathcal{H}}^2 + \left\lVert\frac{v^+}{{\sin(\theta/2)}}\right\rVert_{L^\infty}\llbracket \eta^+ \rrbracket_\mathcal{H}^2\\
        &\lesssim \llbracket \eta^+ \rrbracket_{\mathcal{H}}^3 + \llbracket \eta^+ \rrbracket_{\mathcal{H}}\llbracket \eta^- \rrbracket_{\mathcal{H}}^2,
    \end{align*}
    where we have applied \ref{etaLinfty}, \ref{vLinfty}, and integration by parts. To estimate $I^+_2$
    \begin{align*}
        |I^+_2| &= \frac{|\eta^+_{c,0}|}{2\pi}\left|\int_{\mathbb{T}} \frac{\cos(\theta/2)}{{\sin(\theta/2)}} \left(\eta^+v^+_{\theta\theta} - v^+\eta^+_{\theta\theta} - \eta^-v^-_{\theta\theta} + \eta^-_{\theta\theta}v^- \right)d\theta \right|\\
        &\lesssim |\eta^+_{c,0}| \lVert v^+_{\theta\theta }\rVert_{L^2}\llbracket \eta^+ \rrbracket_\mathcal{H} + |\eta^+_{c,0}| \lVert{v^-_{\theta\theta}}\rVert_{L^2}\llbracket \eta^- \rrbracket_\mathcal{H}\\
        &\quad + |\eta^+_{c,0}|\left|\int_{\mathbb{T}} \frac{\cos\left(\theta/2\right)}{{\sin(\theta/2)}} \left(- v^+\eta^+_{\theta\theta} + \eta^-_{\theta\theta}v^- \right)d\theta \right|\\
        &\lesssim |\eta^+_{c,0}|\left(\llbracket \eta^+ \rrbracket_\mathcal{H}^2 + \llbracket \eta^- \rrbracket_\mathcal{H}^2\right) + |\eta^+_{c,0}| \left(\int_{\mathbb{T}}\frac{|\eta^+_\theta v^+| + |\eta^-_\theta v^-|}{\left|{\sin(\theta/2)}\right|^2} + \frac{|\eta^+_\theta v^+_\theta| + |\eta^-_\theta v^-_\theta|}{\left|{\sin(\theta/2)}\right|} d\theta\right)\\
        &\lesssim |\eta^+_{c,0}|\left(\llbracket \eta^+ \rrbracket_\mathcal{H}^2 + \llbracket \eta^- \rrbracket_\mathcal{H}^2\right),
    \end{align*}
    where we have used integration by parts and similar techniques as in $|I^+_1|$. For $\eta^-\in\mathcal{H}$, we have
    \begin{align*}
        \frac{1}{2}\frac{d}{dt}\llbracket \eta^- \rrbracket_{\mathcal{H}_0}^2 = \langle \eta^-, L^-\eta^-\rangle_{\mathcal{H}_0} + \ep\langle \{\eta^-, v^+\}, \eta^-\rangle_{\mathcal{H}_0} - \ep\langle\{\eta^+, v^-\}, \eta^-\rangle_{\mathcal{H}_0},
    \end{align*}
    we also write the product $\langle N_2, \eta^-\rangle_{\mathcal{H}_0}$ in three terms
    \begin{align*}
        I^-_1 &:= \frac{1}{4\pi}\int_{\mathbb{T}}\frac{1}{\left|{\sin(\theta/2)}\right|^2}\left(\eta^-v^+_{\theta\theta}\eta^-_\theta - \eta^-_{\theta\theta}v^+\eta^-_\theta - \eta^+v^-_{\theta\theta}\eta^-_\theta \right)d\theta,\\
        I^-_2 &:= \frac{1}{4\pi} \int_{\mathbb{T}}\frac{\eta^-_{c,0}\sin(\theta)}{\left|{\sin(\theta/2)}\right|^2}\left(\eta^-v^+_{\theta\theta} - \eta^-_{\theta\theta}v^+ - \eta^+v^-_{\theta\theta} + \eta^+_{\theta\theta}v^- \right)d\theta,\\
        I^-_3 &:= \frac{1}{4\pi}\int_{\mathbb{T}}\frac{1}{\left|{\sin(\theta/2)}\right|^2} \left(\eta^+_{\theta\theta}v^-\eta^-_\theta \right)d\theta,
    \end{align*}
    $|I^-_1|$ and $|I^-_2|$ are estimated as in the case of $|I^+_1|$ and $|I^+_2|$, where we have
    \begin{align*}
        |I^-_1| \lesssim  \llbracket \eta^+ \rrbracket_\mathcal{H}\llbracket \eta^- \rrbracket_\mathcal{H}^2\quad\text{and}\quad |I^-_2|\lesssim |\eta^-_{c,0}|\llbracket \eta^+ \rrbracket_\mathcal{H}\llbracket \eta^- \rrbracket_\mathcal{H}.
    \end{align*}
    Now we estimate $I^+_3$ and $I^-_3$ jointly with
    \begin{align*}
        I^+_3+I^-_3 &= \frac{1}{4\pi}\int_{\mathbb{T}}\frac{1}{\left|{\sin(\theta/2)}\right|^2} \left(\eta^-_{\theta\theta}v^-\eta^+_\theta + \eta^+_{\theta\theta}v^-\eta^-_\theta\right)d\theta \\
        & = \frac{1}{4\pi}\int_{\mathbb{T}}\frac{1}{\left|{\sin(\theta/2)}\right|^2} v^-\partial_\theta (\eta^+_\theta\eta^+_\theta) d\theta\\
        &\lesssim \left|\int_{\mathbb{T}}\frac{v^-_\theta\eta^+_\theta\eta^-_\theta}{\left|{\sin(\theta/2)}\right|^2}d\theta\right|
        + \left|\int_{\mathbb{T}} \frac{\cos\left(\theta/2\right)}{{\sin(\theta/2)}^3} \left(v^-\eta^+_\theta\eta^-_\theta\right) d\theta\right|\\
        &\lesssim \llbracket \eta^+ \rrbracket_\mathcal{H}\llbracket \eta^- \rrbracket^2_\mathcal{H}.
    \end{align*}
    Hence we have shown, for some constant $C>0$
    \begin{align*}
        \frac{1}{2}\frac{d}{dt}\left(\llbracket \eta^+ \rrbracket_{\mathcal{H}_0}^2 + \llbracket \eta^- \rrbracket_{\mathcal{H}_0}^2\right)
        \leq& -\frac{3}{8} \llbracket \eta^+ \rrbracket_{\mathcal{H}_0}^2 - \frac{1}{2} \llbracket \eta^- \rrbracket_{\mathcal{H}_0}^2
        + \ep C\left(\llbracket \eta^+ \rrbracket_\mathcal{H}^3 +  \llbracket \eta^+ \rrbracket_\mathcal{H}\llbracket \eta^- \rrbracket^2_\mathcal{H}\right)\\ &+ \ep C\left(|\eta^+_{c,0}|(\llbracket \eta^+ \rrbracket_\mathcal{H}^2 + \llbracket \eta^- \rrbracket_\mathcal{H}^2) + |\eta^-_{c,0}|\llbracket \eta^+ \rrbracket_\mathcal{H}\llbracket \eta^- \rrbracket_\mathcal{H}\right),
    \end{align*}
    then using \ref{conserve}, we have
    \begin{equation}
        \begin{aligned}
            \frac{1}{2}\frac{d}{dt}\left(\llbracket \eta^+ \rrbracket_{\mathcal{H}_0}^2 + \llbracket \eta^- \rrbracket_{\mathcal{H}_0}^2\right)
            \leq& -\frac{3}{8} \left(\llbracket \eta^+ \rrbracket_{\mathcal{H}_0}^2 + \llbracket \eta^- \rrbracket_{\mathcal{H}_0}^2\right)
            + \ep C\llbracket \eta^+ \rrbracket_{\mathcal{H}_0}^3\\
            &+ \ep C\llbracket \eta^+ \rrbracket_{\mathcal{H}_0} (|\eta^-_{c,0}|^2 + \llbracket \eta^- \rrbracket_{\mathcal{H}_0}^2),
        \end{aligned}
        \label{nonlinearenergyq=0}
    \end{equation}
    \subsection{Nonlinear analysis of \texorpdfstring{$\eta^-_{c,0}$}{eta-{c,0}}}
    From the linear analysis, we see that
    \begin{equation}
        \frac{d}{dt}\eta^-_{c,0}(t) = -2\eta^-_{c,0}(t) + F_1(t) + \ep F_2(t),
        \label{nonlinearZ2q=0}
    \end{equation}
    where $F_2$ is the nonlinear part and $F_1$ is the linear contribution from $e_{s,k}$'s and $e_{c,k}$'s given by
    \begin{align*}
        F_1(t) := \sum\nolimits_{k\geq 1}\tfrac{k^2+3k+1}{k^2(k+1)^2} \eta^-_{c,k}(t).
    \end{align*}
    Then we can bound the norm of $F_1$ by
    \begin{align*}
        |F_1| \leq \left(\sum\nolimits_{k\geq1}\left(\tfrac{k^2+3k+1}{k^2(k+1)^2}\right)^2 \right)^\frac{1}{2} \left(\sum\nolimits_{k\geq1} |\eta^-_{c,k}|^2 \right)^\frac{1}{2} \leq 2\llbracket \eta^- \rrbracket_{\mathcal{H}_0},
    \end{align*}
    To bound $F_2$, note that $Z_2$ is a one dimensional subspace spanned by $e_{c,0} = \cos(\theta) - 1$ and $F_2$ is coefficient of $e_{c,0}$ in $\mathbb{P}_{Z_2}N_2$, then
    \begin{align*}
        F_2 = (N_2)_{c,0} = \langle N_2, e_{c,0}\rangle_\mathcal{H} = \langle \{\eta^-, v^+\} - \{\eta^+, v^-\}, \cos(\theta)-1\rangle_\mathcal{H},
    \end{align*}
    then we can reuse the arguments on $|I^-_2|$ from above and get
    \begin{align*}
        |F_2| \lesssim \llbracket \eta^+ \rrbracket_{\mathcal{H}_0}(|\eta^-_{c,0}| + \llbracket \eta^- \rrbracket_{\mathcal{H}_0}).
    \end{align*}
    \subsection{Exponential decay by a bootstrapping type argument in the case \texorpdfstring{$q=0$}{q=0}}
    We first recall the nonlinear bounds on $\llbracket \eta^\pm \rrbracket_{\mathcal{H}_0}$ and $\eta^-_{c,0}$ from \eqref{nonlinearenergyq=0} and \eqref{nonlinearZ2q=0}, for some $C>0$
    \begin{align*}
        \frac{d}{dt} E^2_{\mathcal{H}_0} &\leq -\frac{3}{4}E^2_{\mathcal{H}_0} + \ep C E^3_{\mathcal{H}_0} + \ep C |\eta^-_{c,0}|^2E_{\mathcal{H}_0},\\
        \frac{d}{dt} \eta^-_{c,0} &= -2\eta^-_{c,0} + F_1(t) + \ep F_2(t),
    \end{align*}
    where $E^2_{\mathcal{H}_0} := \llbracket \eta^+ \rrbracket_{\mathcal{H}_0}^2 + \llbracket \eta^- \rrbracket_{\mathcal{H}_0}^2$, $|F_1| \leq 2E_{\mathcal{H}_0}$, and $|F_2|\leq CE_{\mathcal{H}_0}^2 + CE_{\mathcal{H}_0}|\eta^-_{c,0}|$.
    For $H^2(\mathbb{T})$ initial data $\eta^\pm_0\in Y$, a priori, assume that on some time interval the following bounds hold
    \begin{equation}
        E_{\mathcal{H}_0}^2 < 2\Gamma^2 e^{-2\beta t}\quad\text{and}\quad |\eta^-_{c,0}| < 5\Gamma e^{-\beta t}. \tag{$\ast$}\label{bootstrapq=0}
    \end{equation}
    where we pick $0 < \beta < \frac{3}{8}$ and $\Gamma > \max\left\{\llbracket \eta^+_0 \rrbracket_{\mathcal{H}_0}, \llbracket \eta^-_0 \rrbracket_{\mathcal{H}_0}, \llbracket \eta^-_0 \rrbracket_{Z_2}\right\}$. In particular, let
    \begin{align*}
        T := \sup\{t \in [0,\infty) : \text{ condition \ref{bootstrapq=0} holds on the interval }[0,t)\}.
    \end{align*}
    Then by the local existence theorem of 1-D MHD \ref{existence} and by the choice of $\Gamma$, we see that $T>0$. We will also make the following choice for $\ep$ for which the reason will become clear later
    \begin{align*}
         \ep<\min\left\{\tfrac{3-8\beta}{54\sqrt{2}C\Gamma},\quad \tfrac{2-2\beta}{(2+5\sqrt{2})C\Gamma} \right\}.
    \end{align*}
    Here for the sake of contradiction, assume $T<\infty$. By the continuation criterion \ref{continuation} and the a priori bound, we see that the solution of $\eta^\pm$ exists beyond an interval containing $T$. Then by continuity of $\llbracket \eta^\pm \rrbracket_{\mathcal{H}_0}$ and $\eta^-_{c,0}$, at time $T$, we either have $E^2_{\mathcal{H}_0} = 2\Gamma^2 e^{-2\beta T}$ or $|\eta^-_{c,0}| = 5\Gamma e^{-\beta T}$.
    \par
    Suppose the former is true, considering the time derivative of $E^2_{\mathcal{H}_0}$, then
    \begin{align*}
        \left.\frac{d}{dt} E^2_{\mathcal{H}_0}\right\vert_{t=T} \geq \left.\frac{d}{dt}\left(2\Gamma^2 e^{-2\beta t}\right) \right\vert_{t=T} = -4\beta\Gamma^2 e^{-2\beta T},
    \end{align*}
    but by the a priori bound and by continuity, we have
    \begin{align*}
        \left.\frac{d}{dt} E^2_{\mathcal{H}_0}\right\vert_{t=T} &\leq -\left.\frac{3}{4}E^2_{\mathcal{H}_0} + \ep C E^3_{\mathcal{H}_0} + \ep C |\eta^-_{c,0}|^2E_{\mathcal{H}_0} \right\vert_{t=T}\\
        &\leq -\frac{3}{2}\Gamma^2 e^{-2\beta T} + 27\sqrt{2}\ep C \Gamma^3 e^{-3\beta T},
    \end{align*}
    giving the result $\ep\geq\frac{3-8\beta}{54\sqrt{2}C\Gamma}$, which contradicts our choice of $\ep$.
    \par
    Now to show the latter cannot be true, apply Duhamel's principle on $\eta^-_{c,0}$,
    \begin{align*}
        \eta^-_{c,0}(t) &= e^{-2t}\eta^-_{c,0}(0) + \int_0^t e^{-2(t-s)} \left(F_1(s) + \ep F_2(s)\right) ds,
    \end{align*}
    which is bounded by
    \begin{align*}
        |\eta^-_{c,0}|&\leq e^{-2t}\llbracket \eta^-_0 \rrbracket_{Z_2} + e^{-2 t}\int_0^t e^{2s}\left(2\sqrt{2}\Gamma e^{-\beta s} + 2\ep C\Gamma^2e^{-2\beta s} + 5\sqrt{2}\ep C\Gamma^2 e^{-2\beta s}\right) ds\\
        &\leq \Gamma e^{-\beta t} + \tfrac{16\sqrt{2}}{13}\Gamma e^{-\beta t} + \ep \tfrac{C\Gamma(2+5\sqrt{2})}{2-2\beta}\Gamma e^{-\beta t}.
    \end{align*}
    The bound above yields $|\eta^-_{c,0}|\leq 4\Gamma e^{-\beta t}$ with our choice of $\ep$ for all $t\in [0,T)$, then by continuity, $|\eta^-_{c,0}|$ cannot reach the value of $5\Gamma e^{-\beta T}$ at time $T$. This establishes theorem \ref{thm1}.
    \section{Analysis of the operator \texorpdfstring{$Q$}{Q}}
    In this section, we are going to consider the operator $Q$ defined in \eqref{definition} and the effect of $Q$ on the linear part of equation \eqref{Leta-}. Recall that $Q: L^2(\mathbb{T}) \rightarrow L^2(\mathbb{T})$ is defined as
    \begin{align*}
        Qf = \cos(\theta) f + \sin(\theta) Hf.
    \end{align*}
    \begin{remark}
        It might be interesting to comment that on the Fourier basis $e^{ik\theta}$'s for nonzero $k$'s, $Q$ is a shift operator, and $Q$ is bounded on Sobolev spaces $H^s(\mathbb{T})$ for all $s$.
        \begin{equation*}
            Qe^{ik\theta} = e^{i(k-\operatorname{sgn}(k))\theta}\quad\text{for}\quad k\neq 0.
        \end{equation*}
    \end{remark}
    To see that $Q$ is bounded on $Y = Z_1 \oplus \mathcal{H}$, note that $Q(Z_1) = \{0\}$, so we only need to show that $Q$ is bounded on $\mathcal{H}$. Let $f\in \mathcal{H}$ and consider the decomposition of $Qf$ as the following
    \begin{align*}
        Qf= \sin(\theta)Hf(0) + \cos(\theta)f + \sin(\theta)(Hf(\theta) - Hf(0)),
    \end{align*}
    then the first summand is in $Z_1$, and by \ref{embedding}, we see that $|Hf(0)| \lesssim \llbracket f \rrbracket_{\mathcal{H}}$ so it is bounded in $Z_1$. For the remaining parts, we show that it is bounded in $\mathcal{H}$
    \begin{align*}
        \int_\mathbb{T} \left|\frac{\partial_\theta \left(\mathbb{P}_\mathcal{H} Qf\right)}{{\sin(\theta/2)}}\right|^2 d\theta &\lesssim \int_\mathbb{T} \left|\frac{\cos(\theta)f_\theta - \sin(\theta)f}{{\sin(\theta/2)}}\right|^2 + \left|\frac{\sin(\theta)Hf_\theta}{{\sin(\theta/2)}}\right|^2\\&\quad + \left|\frac{\cos(\theta)(Hf(\theta) - Hf(0))}{{\sin(\theta/2)}}\right|^2 d\theta,
    \end{align*}
    the first two summands are obviously bounded by $\llbracket f \rrbracket_\mathcal{H}^2$, while for the third term, observe that
    \begin{equation}
        \begin{aligned}
            Hf(\theta) - Hf(0) &= \frac{1}{2\pi}\mathrm{p.v.} \int_{-\pi}^{+\pi} \left( \cot\left(\frac{\theta - \vartheta}{2}\right) + \cot\left(\frac{\vartheta}{2}\right)\right)f(\vartheta)d\vartheta\\
            &= \frac{1}{2\pi}\mathrm{p.v.}\int_{-\pi}^{+\pi} \left( \frac{{\sin(\theta/2)}}{\sin\left(\vartheta/2\right)} \csc\left(\frac{\theta-\vartheta}{2}\right) \right) f(\vartheta)d\vartheta\\
            &= \sin\left(\frac{\theta}{2}\right)\frac{1}{2\pi}\mathrm{p.v.}\int_{-\pi}^{+\pi} \csc\left(\frac{\theta-\vartheta}{2}\right) \frac{f(\vartheta)}{\sin\left(\vartheta/2 \right)} d\vartheta,
        \end{aligned}
        \label{hilbert}
    \end{equation}
    then by \ref{etaLinfty}, the $L^2(\mathbb{T})$ norm of the integrand $f(\vartheta)\sin(\vartheta/2)^{-1}$ is bounded by $\llbracket f \rrbracket_\mathcal{H}$, then consider the singular integral operator defined on $\mathbb{T}$ by
        \begin{align*}
            f(\theta) \mapsto \frac{1}{2\pi}\mathrm{p.v.}\int\csc\left(\frac{\theta - \vartheta}{2}\right) f(\vartheta)d\vartheta,
        \end{align*}
    by comparing this operator to the Hilbert transformation, it is apparent that their difference defines a bounded operator from $L^2(\mathbb{T})$ to $L^2(\mathbb{T})$. Hence this concludes our claim.
    \subsection{Behavior of \texorpdfstring{$Q$}{Q} on the basis functions of \texorpdfstring{$Y$}{Y}}
    To obtain a precise bound on the norm of $Q$ on $Y$, in this section, we will examine its behavior on $e_{s,k}$'s and $e_{c,k}$'s. We have already observed that $Qe_{s,0} = 0$, and for $e_{c,0} = \cos(\theta) - 1$, we have $
    Qe_{c,0} = Q(\cos(\theta)-1) = 1-\cos(\theta) = -e_{c,0}$. Hence $Z_1\oplus Z_2$ remains invariant under $L^-$ even in the case $q>0$. Now for $k\geq 2$, on the odd basis functions
    \begin{align*}
        Qe_{s,k} &= Q \left(\frac{\sin((k+1)\theta)}{k+1} - \frac{\sin(k\theta)}{k}\right) = \frac{\sin(k\theta)}{k+1} - \frac{\sin((k-1)\theta)}{k}\\
        &=  \frac{k}{k+1} \left( \frac{\sin(k\theta)}{k} - \frac{\sin((k-1)\theta)}{k-1} \right)\\
        &\quad + \frac{1}{k(k+1)}\sum\nolimits_{j = 2}^{k-1} \left(\frac{\sin(j\theta)}{j} - \frac{\sin((j-1)\theta)}{j-1} \right) + \frac{\sin(\theta)}{k(k+1)}\\
        &=  \frac{k}{k+1} e_{s,k-1} + \frac{1}{k(k+1)}\sum\nolimits_{j = 1}^{k-2} e_{s,j} + \frac{1}{k(k+1)}e_{s,0},
    \end{align*}
    note that this result also holds for all $k=1$ by a simple calculation. For the even basis functions
    \begin{align*}
        Qe_{c,k} &= Q \left(\frac{\cos((k+1)\theta) - 1}{k+1} - \frac{\cos(k\theta) - 1}{k}\right)\\
        &= \frac{\cos(k\theta) - \cos(\theta)}{k+1} - \frac{\cos((k-1)\theta) - \cos(\theta)}{k}\\
        &=  \frac{k}{k+1} \left( \frac{\cos(k\theta)}{k} - \frac{\cos((k-1)\theta)}{k-1} \right)\\
        & \quad + \frac{1}{k(k+1)}\sum\nolimits_{j = 2}^{k-1} \left(\frac{\cos(j\theta)}{j} - \frac{\cos((j-1)\theta)}{j-1} \right) + \frac{2\cos(\theta)}{k(k+1)}\\
        &=  \frac{k}{k+1} e_{c,k-1} + \frac{1}{k(k+1)}\sum\nolimits_{j = 1}^{k-2} e_{c,j} + \frac{2}{k(k+1)}e_{c,0},
    \end{align*}
    with the above calculation holds for all $k\geq2$ and the conclusion also holds for $k=1$. Then to bound $Q$ with respect to the seminorm $\llbracket \cdot \rrbracket_{\mathcal{H}_0}$, we have
    \begin{align*}
        \llbracket Qe_{s,k} \rrbracket_{\mathcal{H}_0} = \llbracket Qe_{c,k} \rrbracket_{\mathcal{H}_0} \leq \tfrac{k}{k+1} + (\mathbbm{1}_{k\geq 2})\tfrac{k-2}{k(k+1)} < 1\quad\text{for all}\quad k\geq 1,
    \end{align*}
    hence given any $\eta^-\in Y$, combining the results above, we have
    \begin{equation*}
        \llbracket Q\eta^- \rrbracket_{\mathcal{H}_0} \leq \llbracket \eta^- \rrbracket_{\mathcal{H}_0}.
    \end{equation*}
    Also, for $Q$'s projection onto $Z_1$ and $Z_2$, we have
    \begin{equation}\label{QonZ1Z2}
        \begin{aligned}
            \mathbb{P}_{Z_1}Q\eta^- &= \sum\nolimits_{k\geq1} \tfrac{1}{k(k+1)}\eta_{s,k}^- \sin(\theta)\\
            \mathbb{P}_{Z_2}Q\eta^- &= \left(\sum\nolimits_{k\geq1}\tfrac{2}{k(k+1)}\eta_{c,k}^--\eta^-_{c,0}\right)(\cos(\theta)-1).
        \end{aligned}
    \end{equation}
    \begin{remark}
        Notice that equation \eqref{QonZ1Z2} agrees with our observation on $\mathbb{P}_{Z_1}Q$ in \eqref{etathetapointwise}, i.e., given any function $f\in Y$, we have that $\mathbb{P}_{Z_1}Qf = (\partial_\theta f(0) + Hf(0))\sin(\theta)$, it is easy to see that they are equivalent by a direct computation.
    \end{remark}
    \section{Nonlinear Analysis in the Case of \texorpdfstring{$0<q<\frac{1}{4}$}{0<q<1/4}}
    In this section, we will investigate the exponential stability of the 1-D MHD \eqref{1dMHD} with $0<q<1/4$. Recall the nonlinear perturbed 1-D MHD
    \begin{align}
        \eta^+_t &= L^+\eta^+ + \ep N_1 = L^+\eta^+ + \ep \{\eta^+, v^+\} - \ep\{\eta^-, v^-\},\\
        \eta^-_t &= L^-\eta^- + \ep N_2 = L^-\eta^- + \ep \{\eta^-, v^+\} - \ep\{\eta^+, v^-\}\notag\\
        & \quad - 2q\ep(\eta^-H\eta^+ - \eta^+H\eta^-).
    \end{align}
    Given that the operator norm of $Q$ with respect to $\llbracket \cdot \rrbracket_{\mathcal{H}_0}$ is bounded by $1$ and the quadratic form $\eta^-\mapsto\langle\eta^-, (L-B)\eta^-\rangle_{\mathcal{H}_0}$ has the bound $\langle\eta^-, (L-B)\eta^-\rangle_{\mathcal{H}_0}\leq -\frac{1}{2} \llbracket \eta^- \rrbracket_{\mathcal{H}_0}$, the choice of $0<q<\frac{1}{4}$ is justified to guarantee the decay of $\eta^-$. For convenience, denote the difference
    \begin{equation*}
        \delta := 1 - 4q > 0.
    \end{equation*}
    \subsection{Estimates of \texorpdfstring{$L^\pm$}{L} in \texorpdfstring{$\llbracket \cdot \rrbracket_{\mathcal{H}_0}$}{H0}}\label{whyq<1/4}
    For the linear part, we have the $\llbracket \cdot \rrbracket_{\mathcal{H}_0}$ estimate on $L^-\eta^-$ by combining the bounds on the operators
    \begin{align*}
        \langle \eta^-, L^-\eta^-\rangle_{\mathcal{H}_0} = \langle \eta^-, \left(L - B\right)\eta^-\rangle_{\mathcal{H}_0} + 2q\llbracket Q\eta^- \rrbracket_{\mathcal{H}_0}\llbracket \eta^- \rrbracket_{\mathcal{H}_0} &\leq -\frac{\delta}{2} \llbracket \eta^- \rrbracket_{\mathcal{H}_0}^2,
    \end{align*}
    and the linear estimate on $L^+\eta^+$ in $\llbracket \cdot \rrbracket_{\mathcal{H}_0}$ remains unchanged.
    \subsection{Estimates of \texorpdfstring{$N_1$}{N1} and \texorpdfstring{$N_2$}{N2} in \texorpdfstring{$\llbracket \cdot \rrbracket_{\mathcal{H}_0}$}{H0}}
    Consider the nonlinear part of $\eta^+$ and let $ J_1^+ + J_2^+ = \langle N_1, \eta^+\rangle_{\mathcal{H}_0}$ with $J_1^+$ and $J_2^+$ defined as
    \begin{align*}
        J_1^+ &= \phantom{-}\langle \{\eta^+, v^+\} - \{\mathbb{P}_\mathcal{H}\eta^-, \mathbb{P}_\mathcal{H} v^-\}, \eta^+\rangle_{\mathcal{H}_0},\\
        J_2^+ &= -\langle \{\mathbb{P}_{Z_1}\eta^-, v^-\} + \{\eta^-, \mathbb{P}_{Z_1}v^-\}, \eta^+\rangle_{\mathcal{H}_0},
    \end{align*}
    since $J_1^+$ is the same as the nonlinear product when $q=0$, we will focus on $J_2^+$
    \begin{align*}
        J_2^+ = \langle -\{\eta^-_{s,0}\sin(\theta), v^-\} + \{\eta^-, \eta^-_{s,0}\sin(\theta)\}, \eta^+\rangle_{\mathcal{H}_0}
        = \eta^-_{s,0}\langle (L+B)\eta^-, \eta^+\rangle_{\mathcal{H}_0}.
    \end{align*}
    Now for $\eta^-$, let $J_1^- + J_2^- - J_3^- + J_4^- = \langle N_2, \eta^-\rangle_{\mathcal{H}_0}$, where similar to the above, are each defined as
    \begin{align*}
        J_1^- &= \langle \{\mathbb{P}_{\mathcal{H}}\eta^-, v^+\} - \{\eta^+, \mathbb{P}_{\mathcal{H}} v^-\}, \eta^+\rangle_{\mathcal{H}_0},\\
        J_2^- &= \langle \{\mathbb{P}_{Z_1}\eta^-, v^+\} - \{\eta^+, \mathbb{P}_{Z_1}v^-\}, \eta^+\rangle_{\mathcal{H}_0},\\
        J_3^- &= 2q\langle \eta^-H\eta^+, \eta^-\rangle_{\mathcal{H}_0},\\
        J_4^- &= 2q\langle \eta^+H\eta^-, \eta^-\rangle_{\mathcal{H}_0},
    \end{align*}
    again, $J^-_1$ is similar to the $q=0$ case, for $J_2^-$, we have
    \begin{align*}
        J_2^- = \langle\{\eta^-_{s,0}\sin(\theta), v^+\} + \{\eta^+, \eta^-_{s,0}\sin(\theta)\}, \eta^- \rangle_{\mathcal{H}_0} = \eta^-_{s,0}\langle (L-B)\eta^+, \eta^-\rangle_{\mathcal{H}_0}.
    \end{align*}
    Hence $|J_1^+ + J_1^-| = |I_1^+-I_2^++I^+_3 + I_1^- - I_2^- + I_3^-| \lesssim \llbracket \eta^+ \rrbracket^3_{\mathcal{H}_0} + \llbracket \eta^+ \rrbracket_{\mathcal{H}_0}(|\eta^-_{c,0}|^2 + \llbracket \eta^- \rrbracket_{\mathcal{H}_0}^2)$ from equation \eqref{nonlinearenergyq=0} and for $J_2^+ + J_2^-$ we have
    \begin{align*}
        |J_2^+ + J_2^-| &= |\eta^-_{s,0}|\left|\langle (L+B)\eta^-, \eta^+\rangle_{\mathcal{H}_0} + \langle (L-B) \eta^+, \eta^-\rangle_{\mathcal{H}_0} \right|\\
        &\leq |\eta^-_{s,0}| \left| \langle L(\eta^++\eta^-), \eta^++\eta^-\rangle_{\mathcal{H}_0}\right| + |\eta^-_{s,0}| \left|\langle L\eta^+, \eta^+\rangle_{\mathcal{H}_0}\right|\\
        &\quad + |\eta^-_{s,0}| \left|\langle L\eta^-, \eta^-\rangle_{\mathcal{H}_0}\right| + 2|\eta^-_{s,0}|\cdot\lVert B \rVert\cdot \llbracket \eta^+ \rrbracket_{\mathcal{H}_0} \cdot \llbracket \eta^- \rrbracket_{\mathcal{H}_0}\\
        &\lesssim |\eta^-_{s,0}|\left(\llbracket \eta^+ \rrbracket_{\mathcal{H}_0}^2 + \llbracket \eta^- \rrbracket_{\mathcal{H}_0}^2\right).
    \end{align*}
    Now we will estimate the term $J_3^-$, which can be written as $2q\langle \eta^-H\eta^+, \mathbb{P}_{\mathcal{H}_0}\eta^-\rangle_{\mathcal{H}}$
    \begin{align*}
        |J_3^-| &\leq \frac{q}{2\pi} \int_\mathbb{T} \frac{1}{|{\sin(\theta/2)}|^2} |\partial_\theta\mathbb{P}_\mathcal{H}\left(\eta^- H\eta^+ \right)| \cdot |\partial_\theta\mathbb{P}_{\mathcal{H}_0} \eta^-| d\theta\\
        &= \frac{q}{2\pi}\int_\mathbb{T} \frac{1}{|{\sin(\theta/2)}|^2} \left|\partial_\theta (\mathbb{P}_\mathcal{H}\eta^- H\eta^+)\right| \cdot |\partial_\theta \mathbb{P}_{\mathcal{H}_0} \eta^-| d\theta\\
        & \quad + \frac{q}{2\pi}\int_\mathbb{T} \frac{1}{|{\sin(\theta/2)}|^2} \left|\partial_\theta (\eta^-_{s,0} \sin(\theta)(H\eta^+(\theta) - H\eta^+(0)) )\right| \cdot |\partial_\theta \mathbb{P}_{\mathcal{H}_0} \eta^-| d\theta\\
        &= \frac{q}{2\pi}\int_\mathbb{T} \frac{1}{|{\sin(\theta/2)}|^2} \left|\mathbb{P}_\mathcal{H}\eta^-_\theta H\eta^+ \right| \cdot |\partial_\theta \mathbb{P}_{\mathcal{H}_0} \eta^-| d\theta\\
        & \quad + \frac{q}{2\pi}\int_\mathbb{T} \frac{1}{|{\sin(\theta/2)}|^2} \left|\eta^-_{s,0}(H\eta^+(\theta) - H\eta^+(0))\cos(\theta)\right| \cdot |\partial_\theta \mathbb{P}_{\mathcal{H}_0} \eta^-| d\theta\\
        & \quad + \frac{q}{2\pi}\int_\mathbb{T} \frac{1}{|{\sin(\theta/2)}|^2} \left|\eta^-H\eta^+_\theta\right| \cdot |\partial_\theta \mathbb{P}_{\mathcal{H}_0} \eta^-| d\theta\\
        &\lesssim \lVert H\eta^+ \rVert_{L^\infty}\llbracket \eta^- \rrbracket_{\mathcal{H}}\llbracket \eta^- \rrbracket_{\mathcal{H}_0} + |\eta^-_{s,0}| \left\lVert\frac{\Delta_\theta H\eta^+}{{\sin(\theta/2)}}\right\rVert_{L^2} \llbracket \eta^- \rrbracket_{\mathcal{H}_0}\\
        & \quad + \lVert H\eta^+_\theta \rVert_{L^2} \left\lVert\frac{\eta^-}{{\sin(\theta/2)}} \right\rVert_{L^\infty} \llbracket \eta^- \rrbracket_{\mathcal{H}_0},
    \end{align*}
    where we have used the term $\Delta_\theta H\eta^+$ to denote the finite difference $H\eta^+(\theta) - H\eta^+(0)$. Then for the term $\lVert \sin(\theta/2)^{-1}\Delta_\theta H\eta^+\rVert_{L^2}$, we use the bound in \eqref{hilbert} to see that it is bounded by $\llbracket \eta^+ \rrbracket_{\mathcal{H}}$, on the other hand, $\lVert\sin(\theta/2)^{-1}\eta^-\rVert_{L^\infty} \lesssim |\eta^-_{s,0}| + \llbracket \eta^- \rrbracket_{\mathcal{H}}$. Then we have
    \begin{align*}
        |J^-_3| &\lesssim \lVert H\eta^+ \rVert_{L^\infty}\llbracket \eta^- \rrbracket_{\mathcal{H}}\llbracket \eta^- \rrbracket_{\mathcal{H}_0} + |\eta^-_{s,0}|\llbracket \eta^+ \rrbracket_{\mathcal{H}} \llbracket \eta^- \rrbracket_{\mathcal{H}_0}\\
        &\quad + \lVert H\eta^+_\theta \rVert_{L^2}(|\eta^-_{s,0}| + \llbracket \eta^- \rrbracket_{\mathcal{H}}) \llbracket \eta^- \rrbracket_{\mathcal{H}_0}\\
        &\lesssim \llbracket \eta^+ \rrbracket_{\mathcal{H}_0} \llbracket \eta^- \rrbracket_{\mathcal{H}_0} \left(|\eta^-_{s,0}| + |\eta^-_{c,0}| + \llbracket \eta^- \rrbracket_{\mathcal{H}_0} \right).
    \end{align*}
    And to estimate $J_4^-$ we have likewise
    \begin{align*}
        |J_4^-| & = \frac{q}{2\pi} \int_\mathbb{T} \frac{1}{|{\sin(\theta/2)}|^2} \left|\partial_\theta\mathbb{P}_\mathcal{H}\left(\eta^+H\eta^-\right)\right| \cdot |\partial_\theta\mathbb{P}_{\mathcal{H}_0}\eta^-| d\theta \\
        &= \frac{q}{2\pi}\int_\mathbb{T} |H\eta^-| \cdot \left|\frac{\eta^+_\theta}{{\sin(\theta/2)}}\right| \cdot \left| \frac{\mathbb{P}_{\mathcal{H}_0}\eta^-_\theta}{{\sin(\theta/2)}}\right| + |H\eta^-_{\theta}|\cdot \left|\frac{\eta^+}{{\sin(\theta/2)}}\right| \cdot \left| \frac{\mathbb{P}_{\mathcal{H}_0}\eta^-_\theta}{{\sin(\theta/2)}}\right| d\theta,\\
        &\lesssim \lVert H\eta^- \rVert_{L^\infty} \left\lVert\frac{\eta^+_\theta}{{\sin(\theta/2)}}\right\rVert_{L^2} \llbracket \eta^- \rrbracket_{\mathcal{H}_0} + \lVert{H\eta^-_{\theta}}\rVert_{L^2} \left\lVert\frac{\eta^+}{{\sin(\theta/2)}}\right\rVert_{L^\infty} \llbracket \eta^- \rrbracket_{\mathcal{H}_0}\\
        &\lesssim \llbracket \eta^+ \rrbracket_{\mathcal{H}_0} \llbracket \eta^- \rrbracket_{\mathcal{H}_0} (|\eta^-_{s,0}| + |\eta^-_{c,0}| + \llbracket \eta^- \rrbracket_{\mathcal{H}_0}).
    \end{align*}
    Hence all together, we have, for $E^2_{\mathcal{H}_0} = \llbracket \eta^+ \rrbracket_{\mathcal{H}_0}^2 + \llbracket \eta^- \rrbracket_{\mathcal{H}_0}^2$ and for some $C>0$.
    \begin{align*}
        \frac{d}{dt}E_{\mathcal{H}_0}^2 \leq - \delta E_{\mathcal{H}_0}^2 + \ep C E_{\mathcal{H}_0}^3 + \ep C(|\eta^-_{s,0}| + |\eta^-_{c,0}|) E_{\mathcal{H}_0}^2.
    \end{align*}
    \subsection{Estimates of \texorpdfstring{$N_1$}{N1} and \texorpdfstring{$N_2$}{N2} in \texorpdfstring{$Z_1$}{Z1}}
    We have already argued that $\mathbb{P}_{Z_1}\eta^+= 0$ for all $t\geq0$ and by \eqref{QonZ1Z2} the linear parts $L^-\eta^-$ has
    \begin{align*}
        (L^-\eta^-)_{s,0} = -2q(Q\eta^-)_{s,0} = -2q\sum\nolimits_{k\geq1} \tfrac{1}{k(k+1)}\eta_{s,k}^-.
    \end{align*}
    Now it remain to examine the $N_2$ term in $Z_1$, it is easy to show that the terms $\{\eta^-, v^+\}$ and $\{\eta^+, v^-\}$ are in $\mathcal{H}$ by referring to the decomposition used in $J_1^-$ and $J_2^-$, and $\eta^+H\eta^-$ is in $\mathcal{H}$ by the comment after \ref{etaLinfty}. As for $\eta^-H\eta^+$, using a similar decomposition as in $J^-_3$, we have
    \begin{align*}
        (N_2)_{s,0} = -2q (\eta^-H\eta^+)_{s,0} = -2qH\eta^+(0)\eta^-_{s,0}.
    \end{align*}
    Hence we have for $\eta^-_{s,0}(t)$,
    \begin{equation}
        \frac{d}{dt}\eta^-_{s,0}(t) = - 2q\ep \eta^-_{s,0}(t)\sum\nolimits_{k\geq1} \tfrac{1}{k(k+1)}\eta_{s,k}^+(t) - 2q\sum\nolimits_{k\geq1} \tfrac{1}{k(k+1)}\eta_{s,k}^-(t)
    \end{equation}
    \begin{remark}
        Notice that the above result agrees with equation \eqref{etathetapointwise}.
    \end{remark}
    \subsection{Estimates of \texorpdfstring{$N_1$}{N1} and \texorpdfstring{$N_2$}{N2} in \texorpdfstring{$Z_2$}{Z2}}
    The treatment of $\mathbb{P}_{Z_2}\eta^+$ is identical to the $q=0$ case, for $\mathbb{P}_{Z_2}\eta^-$, from the $q=0$ case and equation \eqref{QonZ1Z2}, we see that the linear part of the PDE gives
    \begin{align*}
        (L^-\eta^-)_{c,0} &= ((L+B)\eta^-)_{c,0} -2q(Q\eta^-)_{c,0}\\
        &= -(2-2q)\eta^-_{c,0} + \sum\nolimits_{k\geq1}\left(\tfrac{k^2+3k+1}{k^2(k+1)^2} - \tfrac{4q}{k(k+1)}\right)\eta^-_{c,k}.
    \end{align*}
    For the nonlinear contribution of $N_2$, we have
    \begin{align*}
        (N_2)_{c,0} = \langle N_2, e_{c,0}\rangle_\mathcal{H} = F_2 - 2q\langle\eta^-H\eta^+ - \eta^+H\eta^-, \cos(\theta) - 1 \rangle_\mathcal{H},
    \end{align*}
    where $F_2$ is the nonliner contribution from the $q=0$ case. For the second summand, we can use arguments similar to the bounds on $|J_3^-|$ and $|J_4^-|$, which gives
    \begin{align*}
        |\langle\eta^-H\eta^+ - \eta^+H\eta^-, \cos(\theta) - 1 \rangle_\mathcal{H}| \lesssim E^2_{\mathcal{H}_0} + E_{\mathcal{H}_0}|\eta^-_{c,0}| + E_{\mathcal{H}_0}|\eta^-_{s,0}|.
    \end{align*}
    \subsection{Global in time stability by a bootstrapping type argument in the case \texorpdfstring{$0<q<\frac{1}{4}$}{0<q<1/4}}
    Here we will complete the proof of theorem \ref{thm2}. Collecting all the terms above, we have
    \begin{align*}
        \frac{d}{dt}E_{\mathcal{H}_0}^2 &\leq - \delta E_{\mathcal{H}_0}^2 + \ep C E_{\mathcal{H}_0}^3 + \ep C |\eta^-_{c,0}|^2E_{\mathcal{H}_0} + \ep C(|\eta^-_{s,0}| + |\eta^-_{c,0}|) E_{\mathcal{H}_0}^2,\\
        \frac{d}{dt}\eta^-_{s,0} &= - \ep qG_1(t) \eta^-_{s,0} - qG_2(t),\\
        \frac{d}{dt}\eta^-_{c,0} &= -\tfrac{3+\delta}{2}\eta^-_{c,0} + G_3(t) + \ep G_4(t),
    \end{align*}
    for some $C>0$, where $G_1$, $G_2$, and $G_3$ have the following definitions and bounds
    \begin{alignat*}{2}
        &G_1 := \sum\nolimits_{k\geq1} \tfrac{2}{k(k+1)}\eta_{s,k}^+,\quad &&|G_1| \leq \left\lVert\tfrac{2}{k(k+1)}\right\rVert_{l^2}\cdot \lVert{\eta^+_{s,k}\mathbbm{1}_{k\geq 1}}\rVert_{l^2}\leq 3E_{\mathcal{H}_0},\\
        &G_2 := \sum\nolimits_{k\geq1} \tfrac{2}{k(k+1)}\eta_{s,k}^-,\quad &&|G_2|\leq \left\lVert\tfrac{2}{k(k+1)}\right\rVert_{l^2}\cdot \lVert{\eta^-_{s,k}\mathbbm{1}_{k\geq 1}}\rVert_{l^2} \leq 3E_{\mathcal{H}_0},\\
        &G_3 := \sum\nolimits_{k\geq1}\left(\tfrac{k^2+3k+1}{k^2(k+1)^2} - \tfrac{4q}{k(k+1)}\right)\eta^-_{c,k},\quad &&|G_3|\leq |F_1| \leq 2E_{\mathcal{H}_0},
    \end{alignat*}
    and $G_4 := (N_2)_{c,0}$ has the bound
    \begin{align*}
        |G_4| \leq C(E^2_{\mathcal{H}_0} + E_{\mathcal{H}_0}|\eta^-_{c,0}| + E_{\mathcal{H}_0}|\eta^-_{s,0}|).
    \end{align*}
    Now we make an a priori assumption for $\eta^\pm_0\in \mathcal{H}$, let $T$ be
    \begin{align*}
        T := \sup\{t \in [0,\infty) : \text{ condition \ref{bootstrapq>0} holds on the interval }[0,t)\},
    \end{align*}
    where condition \ref{bootstrapq>0} is given by
    \begin{equation}
        E_{\mathcal{H}_0}^2 < 2\Gamma^2 e^{-2\beta t}\quad\text{and}\quad |\eta^-_{c,0}| < 5\Gamma e^{-\beta t}, \tag{$\ast\ast$}\label{bootstrapq>0}
    \end{equation}
    where we pick $\Gamma > \max\left\{\llbracket \eta^+_0 \rrbracket_{\mathcal{H}_0}, \llbracket \eta^-_0 \rrbracket_{\mathcal{H}_0}, \llbracket \eta^-_0 \rrbracket_{Z_2}\right\}$, and in principal, we can pick $\beta$ to be any value given $0<\beta<\frac{\delta}{2}$, here for simplicity we will pick $\beta = \frac{\delta}{4}$. By continuity of norm in $Y$ and the choice of $\Gamma$, we see that $T>0$, and for the sake of contradiction, assume $T<\infty$. 
    We will first make an estimate on $|\eta^-_{s,0}|$ on the interval $[0, T)$ based on the a priori assumption. Let $H$ to be the Green's function generated by $-\ep q G_1$, which is given by
    \begin{align*}
        H(t, s) := \exp\left(-\int_s^t \ep qG_1(\tau) d\tau\right),
    \end{align*}
    then using the a priori bound on $G_1$, we have for any $0\leq s\leq t<T$
    \begin{align*}
        |H(t, s)| \leq \exp\left(\ep q\int^\infty_0 3\sqrt{2}\Gamma e^{-\beta t}\right) \leq \exp\left(\tfrac{3\sqrt{2}\ep q\Gamma}{\beta}\right).
    \end{align*}
    Now we apply Duhamel's principle, which gives
    \begin{align*}
        |\eta^-_{s,0}| &\leq \int_0^t |H(t,s)|\cdot|qG_2(s)| ds
        \leq \exp\left(\tfrac{3\sqrt{2}\ep q\Gamma}{\beta}\right)\int_0^t 3\sqrt{2}q\Gamma e^{-\beta s} ds\\
        &\leq \exp\left(\tfrac{3\sqrt{2}\ep q\Gamma}{\beta}\right)\tfrac{3\sqrt{2} q\Gamma}{\beta},
    \end{align*}
    let us denote the term $C_q := \frac{3\sqrt{2}q}{\beta}$, then we have $|\eta^-_{s,0}| \leq e^{\ep C_q\Gamma} C_q\Gamma$. Now we make the choice for $\ep$
    \begin{align*}
         \ep<\min\left\{\tfrac{1-4q}{12\sqrt{2}q\Gamma},\quad \tfrac{1-4q}{C\Gamma(C_1 + 2eC_q)},\quad \tfrac{3}{2( CC_2\Gamma + eCC_3C_q\Gamma )} \right\},
    \end{align*}
    where $C_1$, $C_2$, $C_3$, and $C_4$ are constants in the bootstrapping estimates independent of $\Gamma$ and $q$ that appears later in this argument. Also note that with this choice of $\ep$ we have $|\eta^-_{s,0}|\leq eC_q\Gamma$.
    By the continuation criterion, the solution of $\eta^\pm$ can be extended beyond time $T$. Then by continuity of $\llbracket \eta^\pm \rrbracket_{\mathcal{H}_0}$ and $\eta^-_{c,0}$, we either have $E^2_{\mathcal{H}_0} = 2\Gamma^2 e^{-2\beta T}$ or $|\eta^-_{c,0}| = 5\Gamma e^{-\beta T}$ at time $T$.
    \par
    Suppose the first case is true, then we consider the time derivative of $E^2_{\mathcal{H}_0}$,
    \begin{align*}
        \left.\frac{d}{dt} E^2_{\mathcal{H}_0}\right\vert_{t=T} \geq \left.\frac{d}{dt}\left(2\Gamma^2 e^{-2\beta t}\right) \right\vert_{t=T} = -4\beta\Gamma^2 e^{-2\beta T},
    \end{align*}
    but by the a priori bound and by continuity, there exists $C_1>0$ such that
    \begin{align*}
        \left.\frac{d}{dt} E^2_{\mathcal{H}_0}\right\vert_{t=T} \leq  - 2\delta\Gamma^2 e^{-2\beta T} + \ep C C_1\Gamma^3e^{-3\beta t} + 2e\ep C C_q\Gamma^3e^{-2\beta T},
    \end{align*}
    giving the result $\ep \geq \frac{1-4q}{C\Gamma(C_1 + 2eC_q)}$, which contradicts our choice of $\ep$.
    \par
    Now to show the second case cannot be true, apply Duhamel's principle on $\eta^-_{c,0}$,
    \begin{align*}
        |\eta^-_{c,0}| &\leq e^{-\frac{3+\delta}{2}t}\llbracket \eta^-_0 \rrbracket_{Z_2} + \int_0^t e^{-\frac{3+\delta}{2}(t-s)} \left|G_3(s) + \ep G_4(s)\right| ds\\
        &\leq e^{-\frac{3+\delta}{2}t}\llbracket \eta^-_0 \rrbracket_{Z_2} + e^{-\frac{3+\delta}{2} t}\int_0^t e^{\frac{3+\delta}{2}s}\left(2\sqrt{2} + \ep C\Gamma\left(C_2 e^{-\beta s} + eC_3 C_q\right)\right)\Gamma e^{-\beta s} ds\\
        &\leq \Gamma e^{-\beta t} + \tfrac{4\sqrt{2}}{3} \Gamma e^{-\beta t} + \tfrac{2}{3}\ep\left(C_2 + eC_3C_q\right) C\Gamma^2 e^{-\beta t}.
    \end{align*}
    This yields $|\eta^-_{c,0}|\leq 4\Gamma e^{-\beta t}$ with the $\ep$ we picked for all $t\in [0,T)$, but by continuity, $|\eta^-_{c,0}|$ cannot reach the value of $5\Gamma e^{-\beta T}$ at time $T$. Hence $T=\infty$ and this shows the exponential decay of $E_{\mathcal{H}_0}$ and $|\eta^+_{c,0}|$. Now let $h(t) := \ep\eta^-_{s,0}(t)$, then we see that $h(t)$ remains bounded, also note that $C_q \in \mathcal{O}(q)$ for $q\rightarrow 0$, then we have that $h\in \mathcal{O}(\ep q)$ for $\ep$ and $q$ small. The stability of $h(t)$ for $t\rightarrow \infty$ is obvious since we have established that the generator $G_1$ and the forcing $G_2$ both decay exponentially in time. This concludes our proof of theorem \ref{thm2}.
    \appendix
    \section{Appendix}
    \subsection{Local existence and continuation criterion of the 1-D MHD}\label{existence}
        \begin{theorem}[Local existence]
            Given initial data $\omega^\pm_0 \in H^k(\mathbb{T})$ for $k\geq1$, then there exists a constant $C > 0$, a time $T > 0$, and a unique solution $\omega^\pm \in C([0, T) ; H^k) \cap \mathrm{Lip}([0, T) ; H^{k-1})$ of the Cauchy problem for the 1-D MHD \eqref{1dMHD}. Moreover, the solution $\omega^\pm$ may be uniquely extended to a maximal time interval $[0,T^\ast)$ where either
            \begin{equation}
                T^\ast = \infty\quad \text{or} \quad \limsup\nolimits_{t\uparrow T^\ast}\displaystyle\left(\lVert \omega^+ \rVert_{H^k} + \lVert \omega^- \rVert_{H^k}\right) = \infty.
            \end{equation}
            \begin{proof}
            Here we will outline the proof following ideas used in establishing local existence of the incompressible Euler equations found in texts such as \cite{Bedrossian2022} and \cite{Hou2009}.
                \begin{enumerate}[leftmargin=1.5em]
                    \item
                    Consider the following a priori energy estimate, suppose there exists a solution, then
                    \begin{align*}
                        \frac{1}{2}\frac{d}{dt}\lVert \omega^\pm \rVert_{H^k}^2 &\lesssim \int_{\mathbb{T}} (1+\partial_\theta^k)\left(p\omega^\pm H\omega^\mp + q\omega^\mp H\omega^\pm - au^\mp \omega^\pm_\theta\right) (1+\partial_\theta^k) \omega^\pm d\theta\\
                        &\lesssim \lVert \omega^\pm \rVert_{H^k}^2\lVert \omega^\mp \rVert_{H^k},
                    \end{align*}
                    hence if we let $E_{H^k} = \left(\lVert \omega^+ \rVert_{H^k}^2 + \lVert \omega^- \rVert_{H^k}^2\right)^\frac{1}{2}$, we have for some $C>0$
                    \begin{equation}
                        \frac{d}{dt} E_{H^k}^2 \leq C E^3_{H^k},\quad\text{so}\quad E_{H^k}(t)\leq \frac{E_{H^k}(0)}{1-CE_{H^k}(0)t}. \label{apriorienergy}
                    \end{equation}
                    \item Consider the mollified 1-D MHD defined by
                    \begin{equation}
                        \partial_t\omega^\pm_\ep + a\varphi_\ep\ast\left(u^\mp_\ep \partial_\theta \left( \varphi_\ep \ast \omega^\pm_\ep\right)\right) = p\omega^\pm_\ep H\omega^\mp_\ep + q\omega^\mp_\ep H\omega^\pm_\ep,
                    \end{equation}
                    where $\{\varphi_\ep\}_{\ep>0}$ is any family of $C^\infty_c$ mollifier defined through scaling $\ep$. Then let
                    \begin{align*}
                        F^\pm_\ep(\omega^\pm_\ep) := p\omega^\pm_\ep H\omega^\mp_\ep + q\omega^\mp_\ep H\omega^\pm_\ep - a\varphi_\ep\ast\left(u^\mp_\ep \partial_\theta \left(\varphi_\ep\ast \omega^\pm_\ep\right)\right)
                    \end{align*}
                    and consider the Hilbert space valued ODE system defined by
                    \begin{align*}
                        \frac{d}{dt}\omega^\pm_\ep (t) = F^\pm_\ep(\omega^\pm_\ep(t)),\quad \omega^\pm_\ep(0) = \varphi_\ep \ast \omega^\pm_0.
                    \end{align*}
                    A routine application of the Banach fixed point theorem shows the existence of the mollified family $\left\{\omega^\pm_\ep\right\}_\ep$ that satisfies the initial condition $\varphi_\ep\ast\omega^\pm_0$ in the space $C([0, T_\ep) ; H^k)\cap \mathrm{Lip}([0,T_\ep) ; H^{k-1})$ where the existence times $T_\ep$'s satisfy $T_\ep \propto \ep E_{H_k}(0)^{-1}$.
                    \item 
                    Given the structure of the mollified equation, we observe that the a priori energy estimate \eqref{apriorienergy} applies on the mollified solutions, then there exists a $T>0$ such that the solutions $\omega^\pm_\ep$'s exist on $[0,T]$ for all $\ep>0$ and are uniformly bounded in the space $C([0, T] ; H^k) \cap \mathrm{Lip}([0, T] ; H^{k-1})$, $T$ can be picked to be any value such that $T<CE_{H_k}(0)^{-1}$.
                    \item 
                    By applying Aubin-Lions compactness theorems, we can find a subsequence of the mollified solutions that converges strongly in the space $C([0,T] ; H^{k-1})$, denote the limit as $\omega^\pm$, then by  Helly's theorem of compactness \cite{Lax2002}, the subsequence can be picked to converge weak$^\ast$ to $\omega^\pm$ in $L^\infty([0, T] ; H^k)$ as well as $\mathrm{Lip}([0,T]; H^{k-1})$.
                    \item
                    It is easy to check $\omega^\pm$ solves \eqref{1dMHD} by approximating using $\omega^\pm_\ep$. Uniqueness follows from Gr\"{o}nwall's inequality. To show $\omega^\pm \in C([0,T] ; H^k)$, we can prove that $\omega^\pm$ is weakly continuous in time and $\lVert \omega^\pm \rVert_{H^k}$ is continuous in time by using regularity of $\omega^\pm_\ep$, then weak continuity and continuity of norm concludes norm continuity.
                    \item
                    The statement involving the maximal interval of existence can be shown by iterating the local existence theorem whenever the $H^k(\mathbb{T})$ norm of $\omega^\pm$ remains bounded.
                \end{enumerate}
            \end{proof}
        \end{theorem}
        \begin{theorem}[Continuation criterion]\label{continuation}
            Given initial data $\omega^\pm_0 \in H^k(\mathbb{T})$ for $k\geq1$, then either the maximal existence time $T^\ast = \infty$ or
            \begin{equation}
                \int_0^{T^\ast} \left( \lVert \omega^+ \rVert_{L^\infty} + \lVert \omega^- \rVert_{L^\infty} + \lVert H\omega^+ \rVert_{L^\infty} + \lVert H\omega^- \rVert_{L^\infty} \right) dt = \infty.\label{criterion}
            \end{equation}
            \begin{proof}
                We will briefly comment on the proof idea, which is similar to Theorem 1.2 and Theorem 1.3 in \cite{Dai2023}. We will show this claim by assuming $T^\ast<\infty$ and \eqref{criterion} false, then the solution can be extended beyond $T^\ast$. We first show the $H^1$ energy is finite and exponentially bounded. Let $E_{L^\infty} := \lVert \omega^+ \rVert_{L^\infty} + \lVert \omega^- \rVert_{L^\infty} + \lVert H\omega^+ \rVert_{L^\infty} + \lVert H\omega^- \rVert_{L^\infty}$, then consider the energy estimate
                \begin{align*}
                    \frac{d}{dt}\left(\lVert \omega^+ \rVert_{H^1}^2 + \lVert \omega^- \rVert_{H^1}^2\right) \lesssim E_{L^\infty}(t) \left(\lVert \omega^+ \rVert_{H^1}^2 + \lVert \omega^- \rVert_{H^1}^2\right),
                \end{align*}
                then it follows from Gr\"{o}nwall's inequality that there exist a constant $C>0$ such that
                \begin{equation}
                    \lVert \omega^+ \rVert_{H^1}^2(t) + \lVert \omega^- \rVert_{H^1}^2(t) \leq \exp\left(C\int_0^t E_{L^\infty}(s) ds\right) \left(\lVert \omega^+_0 \rVert_{H^k}^2 + \lVert \omega^-_0 \rVert_{H^k}^2\right),
                \end{equation}
                for all $t\leq T^\ast$. Then we apply an induction argument on $k$, suppose for the $H^{k-1}(\mathbb{T})$ norm
                \begin{equation}
                    \sup\nolimits_{t\in[0, T^\ast]}\displaystyle \left(\lVert \omega^+(t) \rVert^2_{H^{k-1}} + \lVert \omega^-(t) \rVert^2_{H^{k-1}}\right) < \infty
                \end{equation}
                then apply the $H^k$ energy estimate with $E_{H^{k-1}}$ being the $H^{k-1}(\mathbb{T})$ energy of $\omega^\pm$
                \begin{align*}
                    \frac{d}{dt}\left(\lVert \omega^+ \rVert_{H^k}^2 + \lVert \omega^- \rVert_{H^k}^2\right) \lesssim E_{H^{k-1}}(t) \left(\lVert \omega^+ \rVert_{H^k}^2 + \lVert \omega^- \rVert_{H^k}^2\right),
                \end{align*}
                then applying Gr\"{o}nwall's inequality again shows that for some constant $C>0$ and for $t\leq T^\ast$
                \begin{equation}
                    \begin{aligned}
                        \lVert \omega^+ \rVert_{H^k}^2(t) + \lVert \omega^- \rVert_{H^k}^2(t) &\leq \exp\left( C t \left(\sup\nolimits_{s\in[0, T^\ast]}\displaystyle E_{H^{k-1}}(s)\right)\right) \\
                        &\quad\times\left(\lVert \omega^+_0 \rVert_{H^k}^2 + \lVert \omega^-_0 \rVert_{H^k}^2\right),
                    \end{aligned}
                \end{equation}
                Then $\omega^\pm$ can be extended beyond $T^\ast$ in $H^k(\mathbb{T})$, which concludes our proof.
            \end{proof}
        \end{theorem}
    \subsection{Properties of the Hilbert space \texorpdfstring{$\mathcal{H}$}{H}}
    In this section, we will define $Z_0 := \operatorname{span}_{\mathbb{R}}\{1\}$ to be the space of constant functions on $\mathbb{T}$.
    \begin{lemma}\label{embedding}
        We have the embedding $H^2(\mathbb{T})\hookrightarrow Z_0\oplus Z_1\oplus \mathcal{H} \hookrightarrow H^1(\mathbb{T})$.
        \begin{proof}
            It is apparent that $Z_0\oplus Z_1\oplus \mathcal{H} \hookrightarrow H^1(\mathbb{T})$. To show the embedding of $H^2\hookrightarrow Z_0\oplus Z_1\oplus\mathcal{H}$, let $f\in H^2$, then since $H^2(\mathbb{T})\hookrightarrow C^1(\mathbb{T})$, we can decompose $f$ as
            \begin{align*}
                f(\theta) = f(0) + \partial_\theta f(0)\sin(\theta) + g(\theta)\quad\text{where}\quad g(\theta) = f - f(0)- \partial_\theta f(0)\sin(\theta),
            \end{align*}
            then we have $g(0)=0$ and
            \begin{align*}
                \llbracket g \rrbracket_{\mathcal{H}} &= \frac{1}{2\sqrt{\pi}} \left(\int_{\mathbb{T}}\frac{\left|\partial_\theta f(\theta) - \partial_\theta f(0)\cos(\theta)\right|^2}{\left|{\sin(\theta/2)}\right|^2}d\theta\right)^\frac{1}{2}\\
                &\leq \frac{1}{2\sqrt{\pi}} \left(\int_{\mathbb{T}}\left|\frac{\Delta_\theta f_\theta}{{\sin(\theta/2)}}\right|^2d\theta\right)^\frac{1}{2}
                + \frac{1}{2\sqrt{\pi}} |\partial_\theta f(0)|\left(\int_{\mathbb{T}}\left|\frac{1-\cos(\theta)}{{\sin(\theta/2)}}\right|^2d\theta\right)^\frac{1}{2},
            \end{align*}
            where we have used $\Delta_\theta f_\theta$ to denote the finite difference $\partial_\theta f(\theta) - \partial_\theta f(0)$.
            The second summand is bounded by $\lVert f_\theta \rVert_{L^\infty} \lesssim \lVert f \rVert_{H^2}$. For the first summand, we use Hardy's inequality,
            \begin{align*}
                \left( \int_{\mathbb{T}} \left|\frac{\Delta_\theta f_\theta}{{\sin(\theta/2)}}\right|^2d\theta \right)^\frac{1}{2} 
                &\leq \left( \int_\mathbb{T} \frac{1}{|{\sin(\theta/2)}|^2} \left(\int_0^1 |f_{\theta\theta}(\theta\vartheta)| |\theta| d\vartheta\right)^2 d\theta \right)^\frac{1}{2} \\
                &\lesssim \int_0^1 \left( \int_\mathbb{T} \left|f_{\theta\theta}(\theta\vartheta)\right|^2d\theta \right)^\frac{1}{2} d\vartheta \lesssim \lVert{f_{\theta\theta}}\rVert_{L^2},
            \end{align*}
            hence all together we have $\llbracket g \rrbracket_{\mathcal{H}} \lesssim \lVert f \rVert_{H^2}$.
        \end{proof}
    \end{lemma}
    Here, we will also recall two handy results proven in \cite{Lei2020}
        \begin{lemma}\label{etaLinfty}
            For function $f\in \mathcal{H}$, we have the $L^\infty$ estimate $
                \left\lVert\frac{f}{{\sin(\theta/2)}}\right\rVert_{L^\infty}\lesssim \lVert f \rVert_{\mathcal{H}}$.
            \begin{proof}
                By directly integrating
                \begin{align*}
                    \left|\frac{f(\theta)}{{\sin(\theta/2)}}\right| &= \frac{1}{|{\sin(\theta/2)}|}\left|\int^\theta_0 f_\vartheta(\vartheta) d\vartheta\right|\\
                    &\leq \frac{1}{|{\sin(\theta/2)}|}
                    \left(\int^\theta_0 \sin\left(\frac{\vartheta}{2}\right)^2d\vartheta\right)^\frac{1}{2}
                    \lVert f \rVert_{\mathcal{H}}\\
                    &\lesssim \lVert f \rVert_{\mathcal{H}}.
                \end{align*}
                This also shows that $\mathcal{H}$ defines a Banach algebra and a $H^1(\mathbb{T})$-algebra.
            \end{proof}
        \end{lemma}
        \begin{lemma}\label{vLinfty}
            For function $f\in \mathcal{H}$, we have the estimate  $\left\lVert\frac{v(f)}{{\sin(\theta/2)}}\right\rVert_{L^\infty}\lesssim \lVert f \rVert_{\mathcal{H}}$.
            \begin{proof}
                We directly apply the Sobolev embedding theorem, since $v(f)(0)=0$,
                \begin{align*}
                    \left|\frac{v(f)(\theta)}{{\sin(\theta/2)}}\right|\lesssim \lVert \partial_\theta v(f) \rVert_{L^\infty} \lesssim \lVert \partial_\theta^2 v(f) \rVert_{L^2} \lesssim \lVert f \rVert_{\mathcal{H}},
                \end{align*}
                which finishes the claim.
            \end{proof}
        \end{lemma}

\section*{Acknowledgments}
    The author would like to thank Vlad Vicol and Jiajie Chen for introducing them to this problem and related literature as well as offering support and crucial insights.

%%%%%%%%%%%%%%%%%%%%%%%%%%%%%%%%%%%%%%%%%%%%%%%%%%%%%%
%          7. REFERENCES SECTION
%%%%%%%%%%%%%%%%%%%%%%%%%%%%%%%%%%%%%%%%%%%%%%%%%%%%%%

%       READ THIS SECTION CAREFULLY

% Each of the references below MUST be cited in your article above. Do not include references that are not cited in your article.

% Follow the examples below carefully. We strongly suggest that you copy and paste your reference information directly into our examples.

% List all references in alphabetical order according to the first author's last name.

% Verify each URL works correctly and can be accessed properly. Your URL links should be to reputable websites. The command line for a website link begins with: \url{ }

% Do not add MR or DOI numbers to your references. AIMS production staff will add this information.

% Using BibTex is not recommended but can be handled.
\bibliographystyle{plain}

\begin{thebibliography}{99}
    \bibitem[BV22]{Bedrossian2022}
    Jacob Bedrossian and Vlad Vicol.
    \newblock {\em The mathematical analysis of the incompressible {E}uler and
      {N}avier-{S}tokes equations---an introduction}, volume 225 of {\em Graduate
      Studies in Mathematics}.
    \newblock American Mathematical Society, Providence, RI, [2022] \copyright
      2022.
    
    \bibitem[CC10]{Castro2010}
    A.~Castro and D.~C\'{o}rdoba.
    \newblock Infinite energy solutions of the surface quasi-geostrophic equation.
    \newblock {\em Adv. Math.}, 225(4):1820--1829, 2010.
    
    \bibitem[Che23]{Chen2023}
    Jiajie Chen.
    \newblock On the regularity of the {D}e {G}regorio model for the 3{D} {E}uler
      equations.
    \newblock {\em J. Eur. Math. Soc.}, 2023.
    
    \bibitem[CHH21]{Chen2021}
    Jiajie Chen, Thomas~Y. Hou, and De~Huang.
    \newblock On the finite time blowup of the {D}e {G}regorio model for the 3{D}
      {E}uler equations.
    \newblock {\em Comm. Pure Appl. Math.}, 74(6):1282--1350, 2021.
    
    \bibitem[CLM85]{Constantin1985}
    P.~Constantin, P.~D. Lax, and A.~Majda.
    \newblock A simple one-dimensional model for the three-dimensional vorticity
      equation.
    \newblock {\em Comm. Pure Appl. Math.}, 38(6):715--724, 1985.
    
    \bibitem[DG90]{DeGregorio1990}
    Salvatore De~Gregorio.
    \newblock On a one-dimensional model for the three-dimensional vorticity
      equation.
    \newblock {\em J. Statist. Phys.}, 59(5-6):1251--1263, 1990.
    
    \bibitem[DVZ23]{Dai2023}
    Mimi Dai, Bhakti Vyas, and Xiangxiong Zhang.
    \newblock 1d model for the 3d magnetohydrodynamics.
    \newblock {\em Journal of Nonlinear Science}, 33(5):87, Jul 2023.
    
    \bibitem[EJ20]{Elgindi2020}
    Tarek~M. Elgindi and In-Jee Jeong.
    \newblock On the effects of advection and vortex stretching.
    \newblock {\em Arch. Ration. Mech. Anal.}, 235(3):1763--1817, 2020.
    
    \bibitem[HY09]{Hou2009}
    Thomas~Y. Hou and Xinwei Yu.
    \newblock Introduction to the theory of incompressible inviscid flows.
    \newblock In {\em Nonlinear conservation laws, fluid systems and related
      topics}, volume~13 of {\em Ser. Contemp. Appl. Math. CAM}, pages 1--71. World
      Sci. Publishing, Singapore, 2009.
    
    \bibitem[JSS19]{Jia2019}
    H.~Jia, S.~Stewart, and V.~Sverak.
    \newblock On the de gregorio modification of the constantin--lax--majda model.
    \newblock {\em Archive for Rational Mechanics and Analysis}, 231(2):1269--1304,
      Feb 2019.
    
    \bibitem[KL84]{Kato1984}
    Tosio Kato and Chi~Yuen Lai.
    \newblock Nonlinear evolution equations and the {E}uler flow.
    \newblock {\em J. Funct. Anal.}, 56(1):15--28, 1984.
    
    \bibitem[Lax02]{Lax2002}
    Peter~D. Lax.
    \newblock {\em Functional analysis}.
    \newblock Pure and Applied Mathematics (New York). Wiley-Interscience [John
      Wiley \& Sons], New York, 2002.
    
    \bibitem[LLR20]{Lei2020}
    Zhen Lei, Jie Liu, and Xiao Ren.
    \newblock On the constantin--lax--majda model with convection.
    \newblock {\em Communications in Mathematical Physics}, 375(1):765--783, Apr
      2020.
    
    \bibitem[OSW08]{Wunsch2008}
    Hisashi Okamoto, Takashi Sakajo, and Marcus Wunsch.
    \newblock On a generalization of the {C}onstantin-{L}ax-{M}ajda equation.
    \newblock {\em Nonlinearity}, 21(10):2447--2461, 2008.
\end{thebibliography}

\end{document}